\documentclass[aap,nameyear,MSNbibl,dvips]{arximspdf}
\usepackage{mathrsfs}
\usepackage{mathbh}
\makeatletter
\def\bsuffix #1{#1}
\makeatother
\doi{10.1214/09-AAP623}
\volume{20}
\issue{2}
\pubyear{2010}
\firstpage{430}
\lastpage{461}

\makeatletter

\newtheorem{theorem}{Theorem}[section]
\newtheorem{prop}[theorem]{Proposition}
\newtheorem{cor}[theorem]{Corollary}
\newtheorem{lm}[theorem]{Lemma}

\newcommand{\ud}{d}

\makeatother

\begin{document}
\begin{frontmatter}

\title{The role of the central limit theorem in discovering sharp rates
of convergence to equilibrium
for the solution of the Kac equation}
\runtitle{Kac equation: Sharp rates of convergence to equilibrium}

\begin{aug}
\author[A]{\fnms{Emanuele} \snm{Dolera}\ead[label=e1]{emanuele.dolera@unipv.it}} and
\author[A]{\fnms{Eugenio} \snm{Regazzini}\corref{}\thanksref{t1}\ead[label=e2]{eugenio.regazzini@unipv.it}}
\runauthor{E. Dolera and E. Regazzini}
\affiliation{University of Pavia}
\address[A]{Dipartimento di Matematica ``Felice Casorati''\\
Universit\`a Degli Studi di Pavia \\
Via Ferrata 1\\
27100 Pavia\\
Italy\\
\printead{e1}\\
\phantom{E-mail: }\printead*{e2}} 
\end{aug}

\thankstext{t1}{Supported by MUR, Grant 2006/134525.}

\received{\smonth{6} \syear{2008}}

%
\begin{abstract}
In Dolera, Gabetta and Regazzini [\textit{Ann. Appl. Probab.}
\textbf{19} (2009) 186--201] it is proved that the
total variation distance between the solution $f(\cdot, t)$ of Kac's
equation and the Gaussian density $(0, \sigma^2)$ has an upper bound
which goes to zero with an exponential rate equal to $-1/4$ as $t
\rightarrow+\infty$. In the present paper, we determine a lower bound
which decreases exponentially to zero with this same rate, provided
that a suitable symmetrized form of $f_0$ has nonzero fourth cumulant
$\kappa_4$. Moreover, we show that upper bounds like
$\overline{C}_{\delta} e^{-({1/4})t}\rho_{\delta}(t)$ are
valid for
some $\rho_{\delta}$ vanishing at infinity when $\int_{\mathbb{R}}
|v|^{4+\delta} f_0(v) \,\ud v < +\infty$ for some $\delta$ in $[0, 2[$
and $\kappa_4 = 0$. Generalizations of this statement are presented,
together with some remarks about non-Gaussian initial conditions which
yield the insuperable barrier of $-1$ for the rate of convergence.
\end{abstract}

%

\end{frontmatter}

\section{Introduction}\label{sec1}

In order to determine the rates of relaxation to equilibrium in kinetic
theory, Kac derived the following Boltzmann-like equation, commonly
known as \textit{the Kac equation}:
%
%
\begin{eqnarray}\label{eq:kac}
\frac{\partial f}{\partial t}(v,t) &=& \frac{1}{2\pi}
\int_{0}^{2\pi} \int_{\mathbb{R}}[ f(v \cos\theta- w
\sin\theta,t) \nonumber\\
&&\hspace*{53pt}{} \times
f( v \sin\theta+ w \cos\theta,t) \\
&&\hspace*{79pt}{} - f(v,t) \cdot f(w,t) ] \,\ud w \,\ud\theta
\qquad(v \in\mathbb{R}, t > 0)\hspace*{-26pt}\nonumber
\end{eqnarray}
with some specific probability density function $f_0$ as
initial datum. The resulting Cauchy problem admits a unique solution
within the class of all probability density functions on $\mathbb{R}$.
Such a solution provides the probability distribution at any time of
the velocity of a single particle in a chaotic bath of like molecules
\textit{moving on the real line}; see Kac (\citeyear{K56}, \citeyear{K59})
and McKean (\citeyear{M66}). It
is well known that the probability measure $\mu(\cdot, t)$ determined
by $f(\cdot, t)$ converges to a distinguished Gaussian law in the
\textit{variational metric}, namely
%
%
\begin{equation} \label{eq:convergenza}
\ud_{\mathrm{TV}}(\mu(\cdot, t) ; \gamma_{\sigma}) := {\sup_{B \in\mathscr
{B}(\mathbb{R})}} |\mu(B,t) -
\gamma_{\sigma}(B)| \rightarrow0\qquad (t \rightarrow+\infty)
\end{equation}
where $\gamma_{\sigma}$ denotes the Gaussian distribution
with zero mean and variance $\sigma^2$ and, for any metric space $S$,
$\mathscr{B}(S)$ stands for the Borel class on $S$. It should be
recalled that (\ref{eq:convergenza}) holds true if and only if the
initial datum has finite second moment and $\sigma^2$ is the value of
this moment. The proof of the ``if'' part of this assertion is given in
Dolera (\citeyear{D07}) by adapting arguments explained in Carlen and Lu (\citeyear{CL03}),
whereas the proof of the ``only if'' part is contained in Gabetta and
Regazzini (\citeyear{GR08}).

In regard to the speed of approach to equilibrium, it has been proven that
%
%
\begin{equation} \label{eq:DGR}
\ud_{\mathrm{TV}}(\mu(\cdot, t) ; \gamma_{\sigma}) \leq C_{\ast} e^{-({1/4})t}
\qquad (t \geq0)
\end{equation}
holds, with $C_{\ast}$ being some suitable constant depending
only on the behavior of $f_0$, when $f_0$ has finite fourth moment and
%
%
\begin{equation} \label{eq:cue}
\varphi_0(\xi) := \int_{\mathbb{R}} e^{i \xi x} f_0(x) \,\ud x =
o(|\xi|^{-p})\qquad (|\xi|\rightarrow
+ \infty)
\end{equation}
is valid for some $p > 0$; see Dolera, Gabetta and Regazzini
(\citeyear{DGR09}). This work will be refered to as DGR throughout the rest of the
present paper. Inequality (\ref{eq:DGR}) is known as \textit{McKean's
conjecture} and the above statement constitutes the first satisfactory
support of this conjecture. Other bounds with respect to weak metrics
have been given in Gabetta and Regazzini (\citeyear{GR09}).

At the end of Section 2.2 of DGR, the question of whether the upper
bound in (\ref{eq:DGR}) can be improved is posed. To the best of the
authors' knowledge, this problem has not yet been tackled, except for a
hint on page 370 of Carlen, Carvalho and Gabetta (\citeyear{CCG05}). The main
proposition in the present paper states that the answer is in the
affirmative only in the rather peculiar case in which the \textit{fourth
cumulant} of the density $\tilde{f}_0(x) := \{ f_0(x) + f_0(-x) \}/2$
is zero. The term ``fourth cumulant'' of a probability distribution
$\mathsf{Q}$ on $\mathscr{B}(\mathbb{R})$ refers to the quantity
\[
\kappa_4(\mathsf{Q}) := \int_{\mathbb{R}} (x - \overline{Q})^4
\mathsf
{Q}(\ud x) - 3 \biggl( \int_{\mathbb{R}} (x - \overline{Q})^2 \mathsf
{Q}(\ud x) \biggr)^2,
\]
with $\overline{Q} := \int_{\mathbb{R}} x \mathsf{Q}(\ud x)$,
under the assumption that the fourth moment is finite. This cumulant is
zero, for example, when $\mathsf{Q}$ is Gaussian.

In view of this fact, one could comment on the main proposition by
noting that improvements of the rate expressed by (\ref{eq:DGR}) turn
out to be impossible when $f_0$ is  dissimilar to all of the members in the class
of all Gaussian probability density functions. For the sake of completeness, we recall
that, given the Fourier--Stieltjes transform $q$ of $\mathsf{Q}$, the
$r$th cumulant of $\mathsf{Q}$ is defined to be the coefficient of
$(i\xi)^r/r!$ in the Taylor expansion of $\log(q(\xi))$; see, for
example, Sections 3.14--3.15 of Stuart and Ord (\citeyear{SO87}).

As a further remark on the aforementioned proposition, it is worth
noting its resemblance to well-known facts related to the approximation
of the distribution function $\mathsf{F}_n$ of the ``standardized'' sum
of $n$ independent and identically distributed random variables with
finite variance, by the standard Gaussian distribution $\Phi$. Indeed,
in general, $\mathsf{F}_n$ is approximated by $\Phi$, except for terms
of order $1/\sqrt{n}$. However, higher orders of approximation hold
when the skewness and kurtosis of the common distribution of each
summand are zero. Lyapounov (\citeyear{L01}) was the pioneer of these kinds of
problems, followed by Cram\'{e}r (\citeyear{C37}), Esseen (\citeyear{E45}) and others.

The structure of the paper is as follows. Section \ref{sec2} contains the
presentation of the main results. Section \ref{sec3} deals with the basic
preliminary facts which pave the way for proofs of the main results. It
is split into two subsections. The former consists of a brief
description of the probabilistic interpretation, according to which
$\mu
(\cdot, t)$ can be seen as distribution of a random weighted sum of
random variables. The latter is devoted to the analysis of the error
associated with the approximation of the law of certain weighted sums
of independent random variables to the Gaussian distribution. Section \ref{sec4}
contains the proofs of the main results stated in Section \ref{sec2}. Finally,
some purely technical details are deferred to the \hyperref[app]{Appendix}, together
with the proofs of two lemmas formulated in Section \ref{sec3}.

\section{Presentation of the new results}\label{sec2}

In order to present the main results we intend to prove in this paper,
it is worth mentioning the following weak version of Kac's problem
(\ref{eq:kac}) proposed in Bobylev (\citeyear{B84}). Taking the Fourier
transform of
both sides of (\ref{eq:kac}) yields
%
%
\begin{equation}\label{eq:kacfou}
\frac{\partial\varphi}{\partial t}(\xi,t) = \frac{1}{2\pi}
\int_{0}^{2\pi} \varphi(\xi\cos\theta, t) \cdot
\varphi(\xi\sin\theta, t) \,\ud\theta- \varphi(\xi,t)
\end{equation}
with initial datum $\varphi_0(\xi) := \int_{\mathbb{R}} e^{i
\xi x} f_0(x) \,\ud x$. It should be noted that if $\varphi_0$ is the
Fourier--Stieltjes transform of any (not necessarily absolutely
continuous) probability distribution $\mu_0$ on $\mathscr{B}(\mathbb
{R})$, then (\ref{eq:kacfou}) can be thought of as a new problem which
generalizes (\ref{eq:kac}). In any case, (\ref{eq:kacfou}) admits a
unique solution $\varphi(\cdot, t)$, which characterizes---in the form
of a Fourier--Stieltjes transform---a probability distribution $\mu
(\cdot, t)$ which, throughout the paper, will be said to be a solution
of (\ref{eq:kacfou}). Obviously, in problem (\ref{eq:kac}), one has
$\mu
_0(B) := \int_B f_0(v) \,\ud v$ and $\mu(B, t) := \int_B f(v, t) \,
\ud v$
for every $B$ in $\mathscr{B}(\mathbb{R})$.

In order to formulate the new results exhaustively, let $\mathfrak
{m}_r$ and $\overline{\mathfrak{m}}_r$ denote the $r$th moment and the
absolute $r$th moment of $\mu_0$, respectively, and let $\tilde{\mu}_0$
be the \textit{symmetrized form of} $\mu_0$ defined by
%
%
\begin{equation} \label{eq:simmetrica}
\tilde{\mu}_0(B) := \{ \mu_0(B) + \mu_0(-B) \}/2,\qquad
B \in\mathscr{B}(\mathbb{R}),
\end{equation}
where $-B$ denotes the set $\{ x | {-}x \in B\}$.

A precise statement of the fact that the rate $-1/4$ may be the best
possible one is contained in the following theorem.
\begin{theorem} \label{thm:invalicabilita}
Suppose that $\mu_0$ possesses finite fourth moment
$\mathfrak{m}_4$ and that
$\kappa_4(\tilde{\mu}_0) \neq0$. Moreover, let $\sigma^2$
be the value of $\mathfrak{m}_2$.
There then exists a strictly positive constant $C$,
depending only on the behavior of $\mu_0$, for which
%
%
\begin{equation} \label{eq:invalicabilita}
\ud_{\mathrm{TV}}(\mu(\cdot, t) ; \gamma_{\sigma}) \geq C e^{-({1/4})t}
\end{equation}
holds true for every $t \geq0$.
\end{theorem}

The proof of this theorem, deferred to Section \ref{sec4}, also contains a
precise quantification of $C$. Since
\[
{\sup_{B \in\mathscr{B}(\mathbb{R})}} |\mathsf{P}(B) - \mathsf
{Q}(B)| =
\frac{1}{2} \int_{\mathbb{R}}
|p(x) - q(x)| \,\ud x
\]
is valid whenever $\mathsf{P}$ and $\mathsf{Q}$ are
absolutely continuous probability distributions
with densities $p$ and $q$, respectively, as an immediate consequence
of Theorem \ref{thm:invalicabilita},
it follows that
%
%
\begin{equation} \label{eq:invalicabilita2}
\frac{1}{2} \int_{\mathbb{R}}
\biggl| f(v,t) - \frac{1}{\sigma\sqrt{2\pi}} e^{- {v^2}/({2
\sigma^2})} \biggr|\,\ud v \geq C e^{-({1/4})t} \qquad(t \geq0)
\end{equation}
is true for the solution $f(\cdot, t)$ of (\ref{eq:kac}),
provided that the initial datum
$f_0$ yields a probability measure $\mu_0$ with the same properties as
in Theorem \ref{thm:invalicabilita}.
From (\ref{eq:invalicabilita2}), it plainly follows that any inequality
such as
\[
\int_{\mathbb{R}}
\biggl| f(v,t) - \frac{1}{\sigma\sqrt{2\pi}} e^{- {v^2}/({2
\sigma^2})} \biggr|\,\ud v \leq C_{\ast} e^{-({1/4})t} \rho(t)\qquad
(t \geq0)
\]
is \textit{not} valid when $\rho$ vanishes at infinity. This
clarifies why inequality (\ref{eq:DGR})
can be viewed as sharp.

We now analyze the effect of assuming that $\kappa_4(\tilde{\mu}_0)
= 0$.
\begin{theorem} \label{thm:D+delta}
Consider Kac's equation (\ref{eq:kac}) with initial datum
$f_0$ such that $\overline{\mathfrak{m}}_{4+\delta} <
+\infty$ for some $\delta$ in $[0,2[$ and
$\kappa_4(\tilde{\mu}_0) = 0$. Further, let $\varphi_0$, the
Fourier transform of $f_0$, satisfy the usual tail condition
(\ref{eq:cue}) for some strictly positive~$p$.
There then exist a strictly positive constant $\overline
{C}_{\delta} = \overline{C}_{\delta}(f_0; p)$ and a function
$\rho_{\delta} \dvtx[0,+\infty[\, \rightarrow[0,+\infty[$ which
vanishes at infinity, for which
%
%
\begin{equation} \label{eq:D+delta}
\int_{\mathbb{R}}
\biggl| f(v,t) - \frac{1}{\sigma\sqrt{2\pi}} e^{- {v^2}/({2
\sigma^2})} \biggr|\,\ud v \leq\overline{C}_{\delta} e^{-({1/4})t}\rho_{\delta}(t)\qquad (t \geq0) .
\end{equation}
In particular, if $\delta$ belongs to $]0,2[$,
one can take
%
%
\begin{equation} \label{eq:rho}
\rho_{\delta}(t) = \exp\{ (-3/4 + 2\alpha_{4+\delta})t \}
\end{equation}
with $\alpha_s := \frac{1}{2\pi}\int_{0}^{2\pi} |{\sin
\theta}|^s \,\ud\theta$.
\end{theorem}

Useful information for quantifying $\overline{C}_{\delta}$ can be found
in Sections \ref{sec43}, \ref{sec44} and Appendices \ref{app2} and \ref{app4}.

Since even cumulants $\kappa_{2m}$ of the Gaussian distribution $(0,
\sigma^2)$ vanish for \mbox{$m \geq2$} and ${\sup_{\xi\in\mathbb{R}}}
|\varphi
(\xi, t) - \operatorname{Re}\varphi(\xi, t)| \leq2 e^{-t}$,
one is led to think that the approach to equilibrium of $\mu(\cdot, t)$
might become faster when the symmetrized form of the initial datum
gives an increasing number of zero even cumulants.
\begin{theorem} \label{thm:piuottimo}
Consider problem (\ref{eq:kac}) and maintain the same
notation as before for $f_0$, $\mu_0$, $\tilde{\mu}_0$, $\varphi_0$
and $\alpha_s$. Further, assume that there exist an
integer $\chi$ greater than 2 and a number $\delta$
in $[0,2[$ for which:
\begin{longlist}
\item $\int_{\mathbb{R}} |v|^{2\chi+ \delta} f_0(v) \,\ud v <
+\infty$;
\item the cumulants $\kappa_{2m}$ of $\tilde{f}_0$
vanish for $m = 2, \ldots, \chi$;
\item $\varphi_0$ meets (\ref{eq:cue}) for some
strictly positive $p$.
\end{longlist}
There then exists a strictly positive constant
$\overline{C}_{\chi, \delta} = \overline{C}_{\chi, \delta}(f_0; p)$
for which
%
%
\begin{equation} \label{eq:D+k+delta}
\int_{\mathbb{R}}
\biggl| f(v,t) - \frac{1}{\sigma\sqrt{2\pi}} e^{- {v^2}/({2
\sigma^2})} \biggr|\,\ud v \leq\overline{C}_{\chi, \delta} e^{-(1 -
2\alpha_{2\chi+ \delta})t}\qquad (t \geq0)
\end{equation}
holds true.
\end{theorem}

Useful information for quantifying $\overline{C}_{\chi, \delta}$ can be
found in Section \ref{sec44} and Appendix~\ref{app2}.

It should be noted that, except for the centered Gaussian law, the most
common distributions do not share
condition (ii), at least for large values of $\chi$. Therefore,
it is reasonable to believe that Theorem \ref{thm:invalicabilita}
covers the usual applications.

It would be interesting to check when, under suitable conditions for
the initial distribution, the value
$-1$ for the rate of relaxation to equilibrium is actually obtained.
The following propositions resolve this issue, under the additional
condition that all moments of $\mu_0$ are finite. It therefore remains
to check whether this moment assumption can actually be recovered from
this high order of relaxation to equilibrium. This problem will be
tackled in a forthcoming work.
\begin{prop} \label{prop:menouno}
If $\mu_0$ possesses moments of every order and the
solution $\mu(\cdot, t)$ of (\ref{eq:kacfou}) satisfies
\[
\ud_{\mathrm{TV}}(\mu(\cdot, t); \gamma_{\sigma}) \leq C e^{-t}
\]
for some strictly positive constant $C$, then
%
%
\begin{equation} \label{eq:classespeciale}
\mu_0(\cdot) = \gamma_{\sigma}(\cdot) + o_{\sigma}(\cdot),
\end{equation}
where $o_{\sigma}$ is a finite signed measure
satisfying $o_{\sigma}(A) = -
o_{\sigma}(-A)$ and $\gamma_{\sigma}(A) + o_{\sigma}(A)
\geq0$
for every Borel subset $A$
of $\mathbb{R}$.
\end{prop}

Observe that the Wild formula [cf. (\ref{eq:Wild}) in Section \ref{sec31}]
implies that\break $\ud_{\mathrm{TV}}(\mu(\cdot, t); \gamma_{\sigma}) =
|o_{\sigma}|
e^{-t}$ when the initial datum is of the type (\ref
{eq:classespeciale}). Therefore, if one assumes there exists some $\rho
\dvtx[0,+\infty[\ \rightarrow[0,+\infty[$ vanishing at infinity so that
$\ud_{\mathrm{TV}}(\mu(\cdot, t); \gamma_{\sigma}) \leq C e^{-t} \rho(t)$, then
the total variation $|o_{\sigma}|$
of $o_{\sigma}$ satisfies $|o_{\sigma}| \leq C \rho(t)$ for all
positive $t$, which is tantamount to asserting
that $o_{\sigma}$ is the null measure. This provides a proof for
the following result.
\begin{cor}
If $\mu_0$ has moments of every order and the solution
$\mu(\cdot, t)$ of (\ref{eq:kacfou}) satisfies
\[
\ud_{\mathrm{TV}}(\mu(\cdot, t); \gamma_{\sigma}) \leq C e^{-t} \rho(t)
\]
for some $\rho$ vanishing at infinity and for
some positive constant $C$, then
$\mu(\cdot, t) = \gamma_{\sigma}(\cdot)$ for every $t \geq0$.
\end{cor}

Thus, if all of the moments of $\mu_0$ are finite, then the value for
the rate of convergence to equilibrium that one cannot
sharpen is just $-1$, unless $\mu_0$ is Gaussian.

\section{Preliminaries}\label{sec3}

To pave the way for the proofs of the main statements, this section
presents some necessary preliminary facts
and results. First, it explains the probabilistic meaning of
Wild's series, originally pointed out in
McKean (\citeyear{M66}). Second, it gives new asymptotic expansions for the
characteristic function of weighted sums
of independent and identically distributed random variables, which
complement analogous statements formulated in, for example, Chapter 8
of Gnedenko and Kolmogorov (\citeyear{GK54}), Chapter 6 of Petrov (\citeyear{P75}) and
Section 3.2 of DGR.

\subsection{McKean's interpretation of Wild's sums}\label{sec31}

Following Wild (\citeyear{W51}), one can express the solution $\varphi(\cdot, t)$
of (\ref{eq:kacfou}) as a time-dependent
mixture of characteristic functions, that is,
%
%
\begin{equation} \label{eq:Wild}
\varphi(\xi, t) = \sum_{n \geq1} e^{-t} (1 - e^{-t})^{n-1}
\hat{q}_{n}(\xi; \varphi_0),
\end{equation}
where
\[
\cases{
\hat{q}_{1}(\xi; \varphi_0) := \varphi_0(\xi), \cr
\displaystyle\hat{q}_{n}(\xi; \varphi_0) = \frac{1}{n-1} \sum_{k=1}^{n-1}
\hat{q}_{k}(\xi; \varphi_0) \star\hat{q}_{n-k}(\xi; \varphi_0)\qquad
(n \geq2)}
\]
and $\star$ denotes the so-called \textit{Wild product} defined by
\[
g_{1}(\xi) \star g_{2}(\xi) := \frac{1}{2\pi} \int_{0}^{2\pi}
g_{1}(\xi\cos \theta) \cdot g_{2}(\xi\sin \theta)
\,\ud\theta.
\]

The Wild series, thanks to a symmetry property of the Wild product,
yields a useful decomposition of $\mu(\cdot, t)$ which we will use
later. Such a decomposition involves the symmetrized form $\tilde{\mu}$
of a probability measure $\mu$ defined by $\tilde{\mu}(B) := [\mu
(B) +
\mu(-B)]/2$ for any $B$ in $\mathscr{B}(\mathbb{R})$. It is well known
that if $\mu^{(s)}(\cdot, t)$ denotes the solution of (\ref{eq:kacfou})
with initial datum $\tilde{\mu}_0$ [see (\ref{eq:simmetrica})], then
one can write
%
%
\begin{equation} \label{eq:restosimmetrico}
\mu(\cdot, t) - \mu^{(s)}(\cdot, t) = o_0(\cdot)e^{-t}
\end{equation}
with $o_0(\cdot) := \mu_0(\cdot) - \tilde{\mu}_0(\cdot)$.

The next description of the probabilistic reinterpretation of (\ref
{eq:Wild}) closely follows Section 3.1 of DGR. Accordingly, we
introduce, using exactly the same notation adopted therein, the
measurable space $(\Omega, \mathscr{F})$ as a product, together with
its coordinate
random elements $\nu$, $\tau$, $\mathbf{\theta} := (\theta_n)_{n
\geq
1}$, $\mathbf{\upsilon} := (\upsilon_n)_{n \geq1}$. We then recall the
definitions of the random elements $\delta_j$, $\pi_j$ given in terms
of \textit{McKean trees} and put $\beta= (\nu, \tau, \mathbf{\theta})$.
Concerning the random variables $\pi_j$, recall the fundamental equality
%
%
\begin{equation} \label{eq:sfera}
\sum_{j=1}^{\nu} \pi^{2}_{j} \equiv1,
\end{equation}
which holds true whenever $\tau$ belongs to $\mathbb{G}(\nu)$.

Now, for some fixed initial datum $\mu_0$ for problem (\ref
{eq:kacfou}), define a family $(\mathsf{P}_t)_{t \geq0}$ of
probability measures on $(\Omega, \mathscr{F})$ according to (12) in
DGR. Next, consider the random variable
%
%
\begin{equation} \label{eq:V}
V = \sum_{j=1}^{\nu} \pi_j \upsilon_j
\end{equation}
and note, via the Wild formula, that
\[
\mu(B, t) = \mathsf{P}_t \{V \in B\}\qquad \bigl(B \in\mathscr
{B}(\mathbb{R}), t \geq0\bigr)
\]
$\mu(\cdot, t)$ being the solution of (\ref{eq:kacfou}) with
$\mu_0$ as initial datum.

Consequently, the random variables $\upsilon_n$ turn out to be
\textit{conditionally independent}, given $\beta$, with respect to each
$\mathsf{P}_t$. Moreover, since $\beta$ and $\mathbf{\upsilon}$ are
independent, one can think of the conditional probability distribution
of $V$ given $\beta$ as the distribution of a weighted sum of
independent random variables. Indeed, for any fixed elementary case
$\overline{\omega}$ in $\Omega$, one can define the random variable
%
%
\begin{equation} \label{eq:Vsegnato}
\overline{V}(\cdot) := \sum_{j=1}^{\nu(\overline{\omega})} \pi
_{j}(\overline{\omega}) \upsilon_j(\cdot)
\end{equation}
on $(\Omega, \mathscr{F})$, for which
%
%
\begin{equation} \label{eq:leggicondizionate}
\mathsf{P}_t \{ V \leq x | \beta\}(\overline{\omega}) = \mathsf
{P}_t \{ \overline{V} \leq x \}\qquad
(x \in\mathbb{R}, t \geq0)
\end{equation}
holds $\mathsf{P}_t$-almost surely in $\overline{\omega}$.
This last equality plays a central role in the rest of the paper
since it allows us to work on a finite sum of independent random
variables using typical tools of the \textit{central
limit problem}. In this context, it is important to examine the
behavior of the moments of the random variable $V$. Their evaluation
essentially depends on sums of powers of the $\pi_j$ via the following
identity proven in Gabetta and Regazzini (\citeyear{GR06}):
%
%
\begin{equation} \label{eq:formulaGR}
\mathsf{E}_{t} \Biggl[ \sum_{j=1}^{\nu} | \pi_j |^{m} \Biggr] =
e^{-(1 - 2\alpha_m)t} ,
\end{equation}
$\alpha_m$ being the same as in Section \ref{sec2}.

\subsection{Some asymptotic expansions for the characteristic
function of weighted sums of independent random variables}\label{sec32}

As in Section 3.2 of DGR, the subject to be investigated here is the
behavior of the characteristic function of weighted sums of independent
and identically distributed random variables.
The expansions given here turn out to be more careful than the
analogous ones contained in the aforementioned work since it is now
assumed that the common probability law of the summands possesses
moments of arbitrarily high order. Cumulants will play a central role
in the analysis of the remainder terms. Finally, the study of the
convergence of weighted sums will provide appropriate conditions to
improve the rate of approach to equilibrium for solutions of
(\ref{eq:kac}).

In the rest of this subsection, $(X_j)_{j \geq1}$ stands for a
sequence of independent and identically distributed real-valued random
variables on some probability space $(E, \mathscr{E}, \mathsf{Q})$ with
common nondegenerate distribution $\zeta$ on $(\mathbb{R}, \mathscr
{B}(\mathbb{R}))$. It is assumed that $\zeta$ is
symmetric [that is, $\zeta(B) = \zeta(-B)$ for every Borel set $B$ of
$\mathbb{R}$] and possesses finite moments up to order $k + \delta$,
where $k = 2\chi$, $\chi$ being some integer greater than 1 and
$\delta
$ being an element of the interval $[0,2[$. Denote the $r$th moment and
the absolute $r$th moment of $\zeta$
by $\mathfrak{m}_r$ and $\overline{\mathfrak{m}}_r$, respectively. Note
that the variance $\sigma^2$ of $\zeta$ coincides with $\mathfrak
{m}_2$. Set $\psi(\xi) :=
\int_{\mathbb{R}} e^{i \xi x} \zeta(\ud x)$, which turns out to be an
even real-valued function, and for every positive integer $n$, define
$\{ c_{1,n},
\ldots, c_{n,n} \}$ to be an array of real constants such that
%
%
\begin{equation} \label{eq:sferadiscreta}
\sum_{j=1}^{n} c_{j,n}^{2} = 1
\end{equation}
holds for every $n$. Now, let $V_n$ be the sum of $Y_{1,n},
\ldots, Y_{n,n}$,
where
\[
Y_{j,n} := \frac{1}{\sigma} c_{j,n} X_{j,n}\qquad (j= 1, \ldots,n)
\]
and let $\psi_n$ be the characteristic function
of $V_n$. Consider the $r$th cumulant $\kappa_r$ and recall that, in
general, it can be defined by
%
%
\begin{equation} \label{eq:cumulanti}
\kappa_r = r! \sum_{(\ast)} (-1)^{s-1} \cdot(s-1)! \cdot
\prod_{l=1}^{r} \frac{1}{k_l !} \biggl( \frac{\mathfrak{m}_l}{l!} \biggr)^{k_l}
\qquad(r = 1, \ldots, k)
\end{equation}
where the symbol $(\ast)$ means that the sum is carried out
over all nonnegative integer solutions
$(k_1, \ldots, k_r)$ of equations
\begin{eqnarray*}
k_1 + 2k_2 + \cdots+ rk_r &=& r, \\
k_1 + k_2 + \cdots+ k_r &=& s
\end{eqnarray*}
with the proviso that $0^0 = 1$. Symmetry of $\zeta$ implies
that existing cumulants of odd order are equal to zero.

From a technical fact proved in the Appendix, Section \ref{app1}, after defining $y_0 := \{
[-6 \sigma^2 + (36 \sigma^4 + 12\mathfrak{m}_4)^{1/2}]/\mathfrak
{m}_4 \}
^{1/2}$, one has $\psi(\xi) \geq1/2$ if $|\xi| \leq y_0$ and
%
%
\begin{equation} \label{eq:logtaylor}
\log\psi(\xi) = \sum_{r=1}^{\chi} (-1)^r \frac{\kappa
_{2r}}{(2r)!} \xi
^{2r} + \xi^k \cdot\epsilon_k(\xi),
\end{equation}
where $\epsilon_k$ is continuous on $[-y_0, y_0]$ and
differentiable on $[-y_0, y_0]\setminus\{0\}$. Moreover, this function
satisfies $\epsilon_k(0) = 0$ and $\lim_{\xi\rightarrow0} \varrho
_{k}(\xi)= 0$, with $\varrho_{k}(\xi) := \xi\cdot\epsilon
_{k}^{\prime}(\xi
)$. Consequently, $M_{0}^{(k)} := {\sup_{\xi\in[-y_0, y_0]}} |\epsilon
_k(\xi)|$ and $M_{1}^{(k)} := {\sup_{\xi\in[-y_0, y_0]}} |\varrho
_{k}(\xi)|$ are two finite constants which depend only on the behavior
of the common probability law $\zeta$.

Now, following the same line of reasoning as in Chapter 6 of Petrov
(\citeyear{P75}), we introduce the quantities
%
%
\begin{equation}
\tilde{\lambda}_{r,n} := \frac{\kappa_{2r}}{\sigma^{2r}} \sum_{j=1}^{n}
c_{j,n}^{2r}\qquad (r = 1, \ldots, \chi)
\end{equation}
and define the polynomials
%
%
\begin{equation} \label{eq:P}
\tilde{P}_{r,n}(\xi) := \sum_{(\ast)} \Biggl( \prod_{m=1}^{r}
\frac{1}{k_{m} !} \biggl( \frac{\tilde{\lambda}_{m+1, n}}{(2m+2)!}
\biggr)^{k_m} \Biggr) (-1)^{r+s} \xi^{2(r+s)}
\end{equation}
for $r = 1, \ldots, \chi-1$. In addition, we introduce another
family of functions $\eta_{k,n}$, which will be used to approximate the
characteristic functions $\psi_n$, defined by
%
%
\begin{equation} \label{eq:etakn}
\eta_{k,n}(\xi) = e^{-\xi^2/2} + \sum_{r=1}^{\chi-1} \tilde
{P}_{r,n}(\xi)
\cdot e^{-\xi^2/2} \qquad(\xi\in\mathbb{R}).
\end{equation}
At this stage, we are in a position to state a couple of
preliminary results that play an important role in the rest of the paper.
\begin{lm} \label{lm:approssimazioni4}
Assume that $\chi= 2$ (i.e., $k = 4$) and
$\delta
= 0$. There then exists a positive constant $C_{4}^{\ast}$,
depending only on the behavior of $\zeta$, such that
%
%
\begin{eqnarray} \label{eq:restopunt4}\hspace*{26pt}
|\psi_{n}(\xi) - \eta_{4,n}(\xi)| &\leq& C_{4}^{\ast} \xi^4 e^{-\xi
^2/2} \Biggl[ \xi^4 \sum_{j=1}^{n} c_{j,n}^4 + \sum_{j=1}^{n} c_{j,n}^4
\biggl| \epsilon_4 \biggl(\frac{c_{j,n} \xi}{\sigma} \biggr) \biggr|\Biggr] ,
\\
\label{eq:restopunt4fast}
|\psi_{n}(\xi) - \eta_{4,n}(\xi)| &\leq& C_{4}^{\ast} \xi^4 (1 +
\xi
^4)e^{-\xi^2/2} \Biggl[ \sum_{j=1}^{n} c_{j,n}^6 + \sum_{j=1}^{n}
c_{j,n}^4 \biggl| \epsilon_4 \biggl(\frac{c_{j,n} \xi}{\sigma} \biggr)
\biggr| \Biggr]
\end{eqnarray}
and
%
%
\begin{eqnarray}\label{eq:restopunt4der}
&&|\psi_{n}^{\prime}(\xi) - \eta_{4,n}^{\prime}(\xi)|
\nonumber\\
&&\qquad\leq C_{4}^{\ast} |\xi|^3 (1 + \xi^6) e^{-\xi^2/2}\\
&&\qquad\quad\hspace*{0pt}{}\times \Biggl[ \sum
_{j=1}^{n} c_{j,n}^6 + \sum_{j=1}^{n} c_{j,n}^4 \biggl( \biggl| \epsilon
_4 \biggl(\frac{c_{j,n} \xi}{\sigma} \biggr) \biggr| + \biggl| \varrho
_4 \biggl(\frac{c_{j,n} \xi}{\sigma} \biggr) \biggr| \biggr) \Biggr]\nonumber
\end{eqnarray}
hold true for every $|\xi| \leq A_{4,n} := \sigma y_0
( \sum_{j=1}^{n} c_{j,n}^{4} )^{-1/4}$.
\end{lm}

In (\ref{eq:restopunt4})--(\ref{eq:restopunt4der}), recall that
$\epsilon_4$ is defined by (\ref{eq:logtaylor}) with $k=4$, and $\rho
_4(\xi) := \xi\epsilon_{4}^{\prime}(\xi)$.

For general $k = 2\chi$ and $\delta$, we have the following result.
\begin{lm} \label{lm:approssimazioni}
If $|\xi| \leq A_{k, \delta, n} := \sigma y_0 ( \sum
_{j=1}^{n} c_{j,n}^{4} )^{-1/(k + \delta)}$, then
%
%
\begin{equation} \label{eq:restopunt}
|\psi_{n}(\xi) - \eta_{k,n}(\xi)| \leq C_{k,\delta}^{\ast}
p_{0,k}(\xi
) |\xi|^{k+\delta} e^{-\xi^2/2} \Biggl( \sum_{j=1}^{n}
|c_{j,n}|^{k+\delta} \Biggr)
\end{equation}
and
%
%
\begin{equation} \label{eq:restoderpunt}
|\psi_{n}^{\prime}(\xi) - \eta_{k,n}^{\prime}(\xi)| \leq
C_{k,\delta}^{\ast}
p_{1,k}(\xi) |\xi|^{k-1+\delta} e^{-\xi^2/2} \Biggl( \sum_{j=1}^{n}
|c_{j,n}|^{k + \delta} \Biggr),
\end{equation}
where $C_{k,\delta}^{\ast}$ is a constant
depending only on the behavior of $\zeta$ and $p_{0,k}(\xi)$,
$p_{1,k}(\xi)$ are polynomials whose coefficients depend only
on $k$.
\end{lm}

The proofs of these lemmata are deferred to Section \ref{app2}, in which one
can also find instructions for the evaluation of $C_{4}^{\ast}$,
$C_{k,\delta}^{\ast}$, $p_{0,k}(\xi)$ and $p_{1,k}(\xi)$. Inequalities
(\ref{eq:restopunt}) and (\ref{eq:restoderpunt}) immediately
entail that
%
%
\begin{equation} \label{eq:restoint}
\biggl(\int_{-A_{k, \delta, n}}^{A_{k, \delta, n}} |\psi_{n}(\xi) - \eta
_{k,n}(\xi)|^2 \,\ud\xi\biggr)^{1/2} \leq C_{k,\delta}^{\ast} a_k
\Biggl( \sum_{j=1}^{n} |c_{j,n}|^{k+\delta} \Biggr)
\end{equation}
and
%
%
\begin{equation} \label{eq:restoderint}
\biggl(\int_{-A_{k,\delta,n}}^{A_{k,\delta,n}} |\psi_{n}^{\prime}(\xi)
- \eta
_{k,n}^{\prime}(\xi)|^2 \,\ud\xi\biggr)^{1/2} \leq C_{k,\delta}^{\ast}
a_k \Biggl( \sum_{j=1}^{n} |c_{j,n}|^{k+\delta} \Biggr),
\end{equation}
where $a_k$ is the maximum between $ (\int_{\mathbb{R}}
\xi^{2k}(1+\xi^2)^2 p_{0,k}^{2}(\xi) e^{-\xi^2} \,\ud\xi)^{1/2}$
and $ (\int_{\mathbb{R}} \xi^{2k-2}(1+\xi^2)^2 p_{1,k}^{2}(\xi)
e^{-\xi^2} \,\ud\xi)^{1/2}$.

\section{Proofs of the main results}\label{sec4}

We first prove Theorem \ref{thm:invalicabilita} and then focus on
Proposition \ref{prop:menouno}. In fact, they rest on similar
arguments. We will then provide proofs for Theorems \ref{thm:D+delta}
and \ref{thm:piuottimo} by adapting methods used in Section 4 of DGR.

Before starting, it is worth introducing some new symbols which will be
used hereafter. First, choose a version of the conditional distribution
function $\mathsf{P}_t \{V \leq x | \beta\}$ and call it
$\mathsf{F}^{\ast}(x)$. In view of (\ref{eq:leggicondizionate}), it
does not depend on $t$. $\mathsf{F}^{\ast}(x)[\overline{\omega}]$ will
indicate dependence of $\mathsf{F}^{\ast}(x)$ on a specific sample
point $\overline{\omega}$ in $\Omega$. The Fourier--Stieltjes transform
of $\mathsf{F}^{\ast}(\cdot)[\overline{\omega}]$ will be
designated by
$\varphi^{\ast}(\cdot)[\overline{\omega}]$. Moreover, an integral over
a measurable subset $S$ of $\Omega$ will often be denoted by $\mathsf
{E}[\cdot; S]$. Symbols $\mathfrak{m}_r$ and
$\overline{\mathfrak{m}}_r$ for $\int x^r \mu_0(\ud v)$ and $\int|x|^r
\mu_0(\ud x)$,
respectively, will continue to be used and $\sigma^2$ will designate
the value of $\mathfrak{m}_2$, while $y_0$ will stand for the quantity
$\{[-6\sigma^2 + (36\sigma^4 + 12\mathfrak{m}_4)^{1/2}]/\mathfrak{m}_4
\}^{1/2}$.

\subsection{Proof of Theorem \protect\ref{thm:invalicabilita}}\label{sec41}

Assume, initially, that $\mu_0$ is symmetric. For simplicity, introduce
the \textit{rescaled solution} $\mu_{\sigma}(\cdot, t)$, defined by
$\mu
_{\sigma}(B, t) := \mu(\sigma B, t)$, where $\sigma B := \{y = \sigma x
| x \in B \}$ for every $B$ in the Borel class of $\mathbb{R}$. By
the homogeneity of the total variation distance, we have $\ud_{\mathrm{TV}}(\mu
(\cdot, t); \gamma_{\sigma}) = \ud_{\mathrm{TV}}(\mu_{\sigma}(\cdot, t);
\gamma
)$, where $\gamma$ is shorthand for the standard normal law $\gamma_1$.
Now, thanks to the elementary inequality
%
%
\begin{equation} \label{eq:VT-CHI}
\ud_{\mathrm{TV}}(\mu_{\sigma}(\cdot, t); \gamma) \geq{\frac{1}{2} \sup
_{\xi
\in\mathbb{R}} }|\varphi(\xi/\sigma, t) - e^{-\xi^2/2}|,
\end{equation}
one can employ the expansions given in Section \ref{sec32}. First,
observe that for any small $\varepsilon$ in $]0, \sigma y_0]$, one has
%
%
\begin{eqnarray} \label{eq:pingu}
&&{\sup_{\xi\in\mathbb{R}}} |\varphi(\xi/\sigma, t) - e^{-\xi^2/2}|
\nonumber\\
&&\qquad\geq
|\varphi(\varepsilon/\sigma, t) - e^{-\varepsilon^2/2}|
\nonumber\\[-8pt]\\[-8pt]
&&\qquad= | \mathsf{E}_t \{ \mathsf{E}_t [ e^{i
\varepsilon V /\sigma} | \beta] - e^{- \varepsilon^2 /2}
\} |\nonumber\\
&&\qquad=
\biggl| \int_{\Omega} \{ \varphi^{\ast}(\varepsilon/\sigma
)[\overline{\omega}] - e^{- \varepsilon^2 /2} \} \mathsf{P}_t(\ud
\overline{\omega}) \biggr| .\nonumber
\end{eqnarray}
Next, after fixing any $\overline{\omega}$ in $\Omega$,
substitute $\nu(\overline{\omega})$ for $n$ and $\pi_j(\overline
{\omega
})$ for $c_{j,n}$ $(j = 1, 2, \ldots, n)$ in Lemma \ref
{lm:approssimazioni4}. This way, $\psi_n(\xi)$ changes into $\varphi
^{\ast}(\xi/\sigma)$ and the restriction that Lemma \ref
{lm:approssimazioni4} imposes on $\varepsilon$ becomes $|\varepsilon|
\leq\sigma y_0$ $( \sum_{j=1}^{\nu(\overline{\omega})} \pi
_{j}^{4}(\overline{\omega}) )^{-1/4}$.
Clearly, this bound holds $\mathsf{P}_t$-almost surely for every $t$,
whenever $\varepsilon$ is not greater then $\sigma y_0$.
 Hence, (\ref
{eq:restopunt4}) can be applied with
\[
\eta_4(\xi)[\overline{\omega}] := e^{-\xi^2/2} + \frac{\kappa_4}{4!
\sigma^4} \Biggl( \sum_{j=1}^{\nu(\overline{\omega})} \pi
_{j}^{4}(\overline{\omega}) \Biggr) \xi^4 e^{-\xi^2/2}
\]
in place of $\eta_{4, n}(\xi)$. If $R_{4}^{\ast}(\xi
)[\overline{\omega}]$ stands for $\varphi^{\ast}(\varepsilon
/\sigma
)[\overline{\omega}] - \eta_4(\xi)[\overline{\omega}]$, then the last
member in (\ref{eq:pingu}) can be written as
%
%
\begin{eqnarray} \label{eq:pingu2}
&&\Biggl| \int_{\Omega} R_{4}^{\ast}(\varepsilon)[\overline{\omega}]
\mathsf{P}_t(\ud\overline{\omega}) + \frac{\kappa_4}{4! \sigma^4}
\varepsilon^4 e^{-\varepsilon^2/2} \int_{\Omega} \Biggl( \sum_{j=1}^{\nu
(\overline{\omega})} \pi_{j}^{4}(\overline{\omega}) \Biggr) \mathsf
{P}_t(\ud\overline{\omega}) \Biggr| \nonumber\\
&&\qquad= \biggl| \int_{\Omega} R_{4}^{\ast}(\varepsilon) \,\ud\mathsf{P}_t
+ \frac{\kappa_4}{4! \sigma^4} \varepsilon^4 e^{-\varepsilon^2/2}
e^{-({1/4})t} \biggr| \\
&&\qquad\geq\biggl| \frac{|\kappa_4|}{4! \sigma^4} \varepsilon^4
e^{-\varepsilon^2/2} e^{-({1/4})t} - \biggl| \int_{\Omega}
R_{4}^{\ast}(\varepsilon) \,\ud\mathsf{P}_t \biggr| \biggr|,\nonumber
\end{eqnarray}
where the equality follows from (\ref{eq:formulaGR}) and the
inequality follows from $ | a + b | \geq| |a| -
|b| |$. Now, the claim is that there exists an $\varepsilon$
independent of $t$ and small enough to have
%
%
\begin{equation} \label{eq:pingu3}
\biggl| \int_{\Omega} R_{4}^{\ast}(\varepsilon) \,\ud\mathsf{P}_t
\biggr| \leq\frac{|\kappa_4|}{2 \cdot4! \sigma^4} \varepsilon^4
e^{-\varepsilon^2/2} e^{-({1/4})t}
\end{equation}
for every nonnegative $t$. To this end, recall the following:
that $\epsilon_4$ [see (\ref{eq:logtaylor})] is a continuous function
depending only on the initial datum $\mu_0$ so that $\epsilon_4(0) =
0$; that $|\kappa_4|$ is strictly positive; that the constant
$C_{4}^{\ast} = C_{4}^{\ast}(\mu_0)$ can never be chosen equal to zero.
The inequality
\[
|\epsilon_4(x)| \leq\frac{|\kappa_4|}{4 \cdot4! \sigma^4
C_{4}^{\ast}}
\]
is surely satisfied for every $x$ belonging to a suitable
nondegenerate interval $[-\overline{x}, \overline{x}]$ included in $[-
y_0, y_0]$. Thus, taking (\ref{eq:restopunt4}) into account, one can write
%
%
\begin{eqnarray}
&&\int_{\Omega} \Biggl[ C_{4}^{\ast} \varepsilon^4 e^{-\varepsilon^2/2}
\sum_{j=1}^{\nu(\overline{\omega})} \pi_{j}^{4}(\overline{\omega})
\Biggl| \epsilon_4 \biggl(\frac{\pi_{j}^{4}(\overline{\omega}) \varepsilon
}{\sigma} \biggr) \Biggr| \Biggr] \mathsf{P}_t(\ud\overline{\omega})
\nonumber\\[-8pt]\\[-8pt]
&&\qquad
\leq\frac{|\kappa_4|}{4 \cdot4! \sigma^4} \varepsilon^4
e^{-\varepsilon^2/2} e^{-({1/4})t}\nonumber
\end{eqnarray}
for every $\varepsilon$ in $]0, \sigma\overline{x}]$ and $t
\geq0$. Moreover,
\[
C_{4}^{\ast} \varepsilon^8 e^{-\varepsilon^2/2} e^{-({1/4})t}
\leq\frac{|\kappa_4|}{4 \cdot4! \sigma^4} \varepsilon^4
e^{-\varepsilon^2/2} e^{-({1/4})t}
\]
is valid for every nonnegative $t$, provided that $\varepsilon
$ is chosen not greater than $\overline{\overline{x}} := ( \frac
{|\kappa_4|}{4 \cdot4! C_{4}^{\ast} \sigma^4} )^{1/4}$. Thus, in
view of (\ref{eq:restopunt4}), (\ref{eq:pingu3}) is satisfied for
$\varepsilon$ in $]0, \min\{ \sigma\overline{x}; \overline
{\overline
{x}} \}]$.

To conclude the proof in the symmetric case, fix $\varepsilon
$ as above in order to have (\ref{eq:pingu3}) and use the following
elementary fact: if $|b| \leq|a|/2$, then $ | |a| - |b| | =
|a| - |b| \geq|a|/2$. Applying this to (\ref{eq:pingu2}), we get
\[
\biggl| \frac{|\kappa_4|}{4! \sigma^4} \varepsilon^4 e^{-\varepsilon
^2/2} e^{-({1/4})t} - \biggl| \int_{\Omega} R_{4}^{\ast
}(\varepsilon) \,\ud\mathsf{P}_t \biggr| \biggr| \geq\frac{|\kappa
_4|}{2 \cdot4! \sigma^4} \varepsilon^4 e^{-\varepsilon^2/2}
e^{-({1/4})t},
\]
which, in view of (\ref{eq:pingu}), provides a lower bound
for $\ud_{\mathrm{TV}}(\mu(\cdot, t); \gamma_{\sigma})$.
When $\mu_0$ is symmetric, the constant $\tilde{C}$, which appears in
Theorem \ref{thm:invalicabilita}, can be taken to be equal to $\frac
{|\kappa_4|}{4 \cdot4! \sigma^4} \varepsilon^4 e^{-\varepsilon
^2/2}$ with
$\varepsilon$ in $]0, \min\{ \sigma\overline{x}; \overline
{\overline
{x}} \}]$.\vspace*{1pt}

When $\mu_0$ is not symmetric, we employ its symmetrized form $\tilde
{\mu}_0$ and recall (\ref{eq:restosimmetrico}) to obtain
\begin{eqnarray*}
\bigl| \mu^{(s)}(B, t) - \gamma_{\sigma}(B) \bigr| &=& | \mu(B, t) - o_0(B)
e^{-t} - \gamma_{\sigma}(B) | \\
&\leq&| \mu(B, t) - \gamma_{\sigma}(B) | + 2e^{-t} \\
&\leq&
\ud_{\mathrm{TV}}(\mu(\cdot, t); \gamma_{\sigma}) + 2e^{-t} \bigl(B \in
\mathscr{B}(\mathbb{R})\bigr),
\end{eqnarray*}
which plainly entails
%
%
\begin{equation} \label{eq:distanzasimm}
\ud_{\mathrm{TV}}\bigl(\mu^{(s)}(\cdot, t); \gamma_{\sigma}\bigr) \leq\ud_{\mathrm{TV}}(\mu
(\cdot, t); \gamma_{\sigma}) + 2e^{-t} .
\end{equation}
From the first part of the proof, one can find a constant
$\tilde{C}(\tilde{\mu}_0) \leq2$ for which
\[
\ud_{\mathrm{TV}}\bigl(\mu^{(s)}(\cdot, t); \gamma_{\sigma}\bigr) \geq\tilde
{C}(\tilde
{\mu}_0) e^{-({1/4})t} .
\]
Hence,
\begin{eqnarray*}
\ud_{\mathrm{TV}}(\mu(\cdot, t); \gamma_{\sigma}) &\geq&\ud_{\mathrm{TV}}\bigl(\mu
^{(s)}(\cdot, t); \gamma_{\sigma}\bigr) - 2e^{-t} \\
&\geq&\tilde{C}(\tilde{\mu}_0) e^{-({1/4})t}- 2e^{-t} \geq
\tfrac{1}{2} \tilde{C}(\tilde{\mu}_0) e^{-({1/4})t}
\end{eqnarray*}
holds, provided that $t \geq\hat{t} := - \log[(\tilde
{C}(\tilde{\mu}_0)/4)^{4/3}]$, where $\hat{t}$ is strictly positive. To
conclude the proof, observe that (\ref{eq:invalicabilita}) is valid,
taking, for example,
\[
\tilde{C} = \tilde{C}(\mu_0) := \min\biggl\{\frac{1}{2} \tilde
{C}(\tilde
{\mu}_0); \inf_{t \in[0, \hat{t}]} \,\ud_{\mathrm{TV}}(\mu(\cdot, t);
\gamma
_{\sigma}) \biggr\} .
\]
Finally, $\inf_{t \in[0, \hat{t}]} \ud_{\mathrm{TV}}(\mu(\cdot, t);
\gamma_{\sigma})$ is strictly positive in view of the existence of the
minimum combined with the uniqueness of the solution of Kac's equation.
This point is clarified in Appendix \ref{app3}.

\subsection{Proof of Proposition \protect\ref{prop:menouno}}\label{sec42}

To prove this proposition under the assumption that all of the moments
of $\mu_0$ are finite, it will suffice to prove that all of the
cumulants $\tilde{\kappa}_{2m}$ of even order of $\tilde{\mu}_0$ are
zero for $m = 2, 3, \ldots.$ Thanks to (\ref{eq:distanzasimm}), the
inequality, which appears in the statement of Proposition \ref
{prop:menouno}, can be rewritten as
%
%
\begin{equation} \label{eq:dinuzzo}
\ud_{\mathrm{TV}}\bigl(\mu^{(s)}(\cdot, t); \gamma_{\sigma}\bigr) \leq(C + 2) e^{-t}.
\end{equation}
In view of this fact, we can assume, without real loss of
generality, that $\mu_0$ is symmetric. Then, supposing that $\kappa
_{2m} = 0$ for $m = 2, \ldots, s-1$ and $\kappa_{2s} \neq0$ for some
integer $s$ greater than 2, we have contradicted (\ref{eq:dinuzzo}).

As in the previous subsection, write
%
%
\begin{eqnarray} \label{eq:fry}
2 \ud_{\mathrm{TV}}(\mu(\cdot, t); \gamma_{\sigma})
&\geq&{\sup_{\xi\in
\mathbb{R}}} |\varphi(\xi/\sigma, t) - e^{-\xi^2/2}|
\nonumber\\[-8pt]\\[-8pt]
&\geq&
\biggl| \int_{\Omega} \{ \varphi^{\ast}(\varepsilon
/\sigma
)[\overline{\omega}] - e^{- \varepsilon^2 /2} \} \mathsf{P}_t(\ud
\overline{\omega}) \biggr|,\nonumber
\end{eqnarray}
where $\varepsilon$ is any positive constant not greater than
$\sigma y_0$. Following the general lines of Section \ref{sec32}, define
\[
\eta_{2s}(\xi)[\overline{\omega}] := e^{-\xi^2/2} + (-1)^s \frac
{\kappa
_{2s}}{(2s)! \sigma^{2s}} \Biggl( \sum_{j=1}^{\nu(\overline{\omega})}
\pi
_{j}^{2s}(\overline{\omega}) \Biggr) \xi^{2s} e^{-\xi^2/2} .
\]
After setting $R_{2s}^{\ast}(\xi)[\overline{\omega}] :=
\varphi^{\ast}(\varepsilon/\sigma)[\overline{\omega}] - \eta
_{2s}(\xi
)[\overline{\omega}]$, the last part of (\ref{eq:fry}) becomes
%
%
\begin{eqnarray}\label{eq:fry2}\qquad
&&\Biggl| \int_{\Omega} R_{2s}^{\ast}(\varepsilon)[\overline{\omega}]
\mathsf{P}_t(\ud\overline{\omega}) + (-1)^s \frac{\kappa
_{2s}}{(2s)! \sigma^{2s}} \varepsilon^{2s} e^{-\varepsilon^2/2} \int
_{\Omega} \Biggl( \sum_{j=1}^{\nu(\overline{\omega})} \pi
_{j}^{2s}(\overline{\omega}) \Biggr) \mathsf{P}_t(\ud\overline{\omega
}) \Biggr| \nonumber\\
&&\qquad= \biggl| \int_{\Omega} R_{2s}^{\ast}(\varepsilon) \,
\ud
\mathsf{P}_t + (-1)^s \frac{\kappa_{2s}}{(2s)! \sigma^{2s}}
\varepsilon^{2s}
e^{-\varepsilon^2/2} e^{-(1 - 2\alpha_{2s})t} \biggr| \\
&&\qquad\geq\biggl| \frac{|\kappa_{2s}|}{(2s)! \sigma^{2s}} \varepsilon^{2s}
e^{-\varepsilon^2/2} e^{-(1 - 2\alpha_{2s})t} - \biggl| \int_{\Omega
} R_{2s}^{\ast}(\varepsilon) \,\ud\mathsf{P}_t \biggr| \biggr| .\nonumber
\end{eqnarray}
Now, if $|\varepsilon| \leq\sigma y_0$, an application of
(\ref{eq:restopunt}), with $k = 2s$ and $\delta= 1$ combined with
(\ref
{eq:formulaGR}), yields
%
%
\begin{eqnarray}\label{eq:fry4}
\biggl| \int_{\Omega} R_{2s}^{\ast}(\varepsilon) \,\ud\mathsf{P}_t\biggr|
&\leq&\int_{\Omega} | R_{2s}^{\ast}(\varepsilon) |
\,\ud\mathsf{P}_t \nonumber\\
&\leq& C_{2s, 1}^{\ast} |\varepsilon|^{2s + 1} \bigl[9\bigl(1 + |\varepsilon
|^{h(s)}\bigr)\bigr] e^{-\varepsilon^2/2} e^{-(1 - 2\alpha_{2s + 1})t}
\\
&\leq& C_{2s, 1}^{\ast} |\varepsilon|^{2s + 1} \bigl[9\bigl(1 + (\sigma
y_0)^{h(s)}\bigr)\bigr] e^{-\varepsilon^2/2} e^{-(1 - 2\alpha_{2s + 1})t}\nonumber
\end{eqnarray}
for every nonnegative $t$. Here, $h(s) := 2s^2 - s$ and the
term $[9(1 + |\varepsilon|^{h(s)})]$ is an upper bound for the
polynomial $p_{0,k}$ in (\ref{eq:restopunt}); see also (\ref{eq:tosca})
in the \hyperref[app]{Appendix}. If $\varepsilon$ satisfies the further restriction
\[
|\varepsilon| \leq\frac{1}{2 C_{2s, 1}^{\ast}} \cdot\frac{1}{9(1 +
(\sigma y_0)^{h(s)})} \cdot\frac{|\kappa_{2s}|}{(2s)! \sigma^{2s}},
\]
then one can rewrite (\ref{eq:fry4}) as
%
%
\begin{equation} \label{eq:fry3}
\biggl| \int_{\Omega} R_{2s}^{\ast}(\varepsilon) \,\ud\mathsf{P}_t
\biggr| \leq\frac{|\kappa_{2s}|}{2 \cdot(2s)! \sigma^{2s}} \varepsilon^{2s}
e^{-\varepsilon^2/2} e^{-(1 - 2\alpha_{2s})t} .
\end{equation}
Hence, inequalities (\ref{eq:fry2}) and (\ref{eq:fry3}) entail
that
\[
\frac{|\kappa_{2s}|}{2 \cdot(2s)! \sigma^{2s}} \varepsilon^{2s}
e^{-\varepsilon^2/2} e^{-(1 - 2\alpha_{2s})t} \leq2 \ud_{\mathrm{TV}}(\mu
(\cdot, t); \gamma_{\sigma}) \leq2 (C + 2) e^{-t}
\]
for every nonnegative $t$, which contradicts the fact that
$(1 - 2\alpha_{2s})$ is strictly smaller than 1. Thus, $\kappa_{2s}$
must vanish, implying that $\mu_0 = \gamma_\sigma$ since $\gamma
_\sigma
$ is uniquely determined by its moments. Finally, if $\mu_0$ is not
symmetric, then $\tilde{\mu}_0 = \gamma_\sigma$.

\subsection{Proof of Theorem \protect\ref{thm:D+delta} when $k +
\delta= 4$}\label{sec43}

We shall closely follow the proof of Theorem 2.1 in DGR. First, let us
assume that the condition
\renewcommand{\theequation}{H}
\begin{equation}\label{hypoH}
f_0 \mbox{ \textit{and}, \textit{consequently}, } f(\cdot,t)
\mbox{ \textit{are even functions}}
\end{equation}
holds. This does not limit the generality of subsequent
reasoning, thanks to (9)--(10) of DGR. Since $\frac{\ud}{\ud v}
\mathsf
{F}^{\ast}(v)$ represents a version of the conditional probability
density function of $V$ given $\beta$, in view of basic properties of
conditional expectation, one has
%
%
\setcounter{equation}{43}
\renewcommand{\theequation}{\arabic{equation}}
\begin{eqnarray}\label{eq:primastima}
&&\int_{\mathbb{R}}
\biggl| f(v,t) - \frac{1}{\sigma\sqrt{2\pi}} e^{- {v^2}/({2
\sigma^2})} \biggr|\,\ud v\nonumber\\
&&\qquad =: \Vert f(v,t) - g_{\sigma}(v) \Vert_{1}
\leq\mathsf{E}_t
\biggl[ \biggl\Vert\frac{\ud}{\ud v} \mathsf{F}^{\ast}(v) - g_{\sigma}(v)
\biggr\Vert_{1}\biggr] \\
&&\qquad=  \mathsf{E}_t \biggl[ \biggl\Vert\frac{\ud}{\ud v}
\mathsf{F}^{\ast}(\sigma v) - g_{1}(v) \biggr\Vert_{1}\biggr],\nonumber
\end{eqnarray}
where $g_{\sigma}(v) \,\ud v = \gamma_{\sigma}(\ud v)$.
Moreover, from Proposition 2.2 of DGR, which can be applied to $f_0$,
thanks to the hypotheses in Theorem \ref{thm:D+delta} and (\ref{hypoH}), there
exist $\alpha$ and $\lambda$ for which
%
%
\begin{equation} \label{eq:magggamma}
|\varphi_0(\xi) |\leq\biggl( \frac{\lambda^2}{\lambda^2
+ \xi^2} \biggr)^{\alpha}
\end{equation}
holds true for every real $\xi$. In particular, one can set
$\alpha= (2 \cdot\lceil2/p \rceil)^{-1}$, $p$ being the same as in
(\ref{eq:cue}) and $\lceil s \rceil$ standing for the least integer not
less than $s$. Define $U \subset\Omega$ by
%
%
\begin{equation} \label{eq:insiemecattivo}
U := \{ \nu\leq\overline{n}\} \cup\Biggl\{ \prod_{j=1}^{\nu}
\pi_j = 0\Biggr\} \cup\Biggl\{ \sum_{j=1}^{\nu} \pi_{j}^{4} \geq
\overline{\delta} \Biggr\}
\end{equation}
with $\overline{n} = 17 \cdot\lceil2/p \rceil$ and
\[
\overline{\delta} = \min\biggl\{ \frac{1}{2^{\overline{n}} \overline
{n}!}; \frac{\sigma^8}{16 y_{0}^{4} \overline{\mathfrak{m}}_{3}^{4}}
\biggr\} \leq\frac{1}{2^{\overline{n}} \overline{n}!} .
\]
Next, check that $U$ belongs to $\mathscr{F}$ and rewrite the last
term in (\ref{eq:primastima}) as
%
%
\begin{equation} \label{eq:secondastima}
\mathsf{E}_t \biggl[ \biggl\Vert\frac{\ud}{\ud v}
\mathsf{F}^{\ast}(\sigma v) - g_{1}(v) \biggr\Vert_{1}; U \biggr] + \mathsf{E}_t
\biggl[ \biggl\Vert\frac{\ud}{\ud v}
\mathsf{F}^{\ast}(\sigma v) - g_{1}(v) \biggr\Vert_{1}; U^{c} \biggr] .
\end{equation}
By the same arguments as the ones used to prove (22) in DGR,
one obtains
\[
\mathsf{P}_t \{ \nu\leq\overline{n} \} \leq\overline{n} e^{-t}
\quad\mbox{and}\quad \mathsf{P}_t \Biggl\{ \prod_{j=1}^{\nu
} \pi_j = 0 \Biggr\} = 0 .
\]
As for the third component of the union in the definition of
$U$, one can combine Markov's (with power 2) and Lyapunov's
inequalities to get
\[
\mathsf{P}_t \Biggl\{ \sum_{j=1}^{\nu} \pi_{j}^{4} \geq\overline{\delta}
\Biggr\} \leq\frac{1}{\overline{\delta}{}^2} \mathsf{E}_t \Biggl[ \Biggl(
\sum_{j=1}^{\nu} \pi_{j}^{4} \Biggr)^2 \Biggr]
\leq\frac{1}{\overline{\delta}{}^2} \mathsf{E}_t \Biggl[ \sum
_{j=1}^{\nu} \pi_{j}^{6} \Biggr] \leq\frac{1}{\overline{\delta}{}^2}
e^{- {(3/8)} t} .
\]
The exponent $3/8$ follows from the application of (\ref
{eq:formulaGR}) with $m = 6$. Now, combining all of the above
computations leads to an estimate for the probability of $U$ under
$\mathsf{P}_t$, that is,
%
%
\begin{equation} \label{eq:probinsiemecattivo}
\mathsf{P}_t (U) \leq[\overline{n} + 1/\overline{\delta}{}^2] e^{-{(3/8)}t}\qquad (t \geq0) .
\end{equation}
Inequality (\ref{eq:probinsiemecattivo}) leads immediately to
the upper bound
%
%
\begin{equation} \label{eq:stimasulbrutto}
\mathsf{E}_t \biggl[ \biggl\Vert\frac{\ud}{\ud v}
\mathsf{F}^{\ast}(\sigma v) - g_{1}(v) \biggr\Vert_{1}; U \biggr] \leq2\mathsf
{P}_t (U) \leq2 [\overline{n} + 1/\overline{\delta}{}^2]
e^{-{(3/8)}t} .
\end{equation}
To control the integral over $U^c$ appearing in (\ref
{eq:secondastima}), we invoke the \textit{Beurling inequality} formulated
in Proposition 4.1 of DGR to obtain
%
%
\begin{eqnarray} \label{eq:beurlingsulbrutto}
&&\mathsf{E}_t \biggl[ \biggl\Vert\frac{\ud}{\ud v}
\mathsf{F}^{\ast}(\sigma v) - g_{1}(v) \biggr\Vert_{1}; U^{c}
\biggr]\nonumber\\[-8pt]\\[-8pt]
&&\qquad \leq\frac
{1}{\sqrt{2}} \mathsf{E}_t \biggl[ \biggl\{ \int_{\mathbb
{R}} |\Delta|^2 \,\ud
\xi+ \int_{\mathbb{R}} |\Delta'|^2 \,\ud\xi\biggr\}^{1/2} ;
U^{c}\biggr],\nonumber
\end{eqnarray}
where\vspace*{1pt} $\Delta:= \varphi^{\ast}(\xi/\sigma) - e^{-
\xi^2/2}$ and $\Delta' := \frac{\ud}{\ud\xi} \Delta$.
Applicability of
this result is justified by the fact that the restriction to $U^c$ of
the conditional characteristic function $\xi\mapsto\varphi^{\ast
}(\xi
) := \int_{\mathbb{R}} e^{i \xi x} \,\ud\mathsf{F}^{\ast}(x)$
belongs to
$H^1(\mathbb{R})$. To see this, note that $\varphi^{\ast}(\xi
)[\overline
{\omega}] = o(|\xi|^{-34})$ is valid for $|\xi| \rightarrow+\infty$
and for $\overline{\omega}$ in $U^c$. Indeed, thanks to conditional
independence and (\ref{eq:magggamma}), one has
\[
|\varphi^{\ast}(\xi)| \leq\prod_{j=1}^{\overline{n}} \biggl( \frac
{\lambda^2}{\lambda^2 + \pi_{j}^2 \xi^2} \biggr)^{\alpha}
\]
and the claimed ``tail behavior'' of $\varphi^{\ast}$ follows
from the definitions of $\overline{n}$ and $\alpha$, together with the
fact that the random numbers $\pi_{j}$ do not vanish on $U^{c}$. To
complete the argument for $H^1(\mathbb{R})$ regularity, use Remark A.2
in Section A.3 of the Appendix of DGR.

Now, the expectation in the right-hand side of (\ref
{eq:beurlingsulbrutto}) is dominated by
%
%
\begin{eqnarray} \label{eq:somme}\qquad
&&\mathsf{E}_t \biggl[ \biggl( \int_{\mathbb{\{|\xi| \leq
\mathrm{A} \}}} |\Delta|^2 \,\ud\xi\biggr)^{1/2} ; U^{c} \biggr] +
\mathsf{E}_t \biggl[ \biggl( \int_{\mathbb{\{|\xi| \geq\mathrm{A} \}}}
|\Delta|^2 \,\ud\xi
\biggr)^{1/2} ; U^{c} \biggr] \nonumber\\[-8pt]\\[-8pt]
&&\qquad{} + \mathsf{E}_t \biggl[ \biggl( \int_{\mathbb{\{|\xi| \leq\mathrm{A}
\}}} |\Delta'|^2 \,\ud\xi\biggr)^{1/2} ; U^{c} \biggr] +
\mathsf{E}_t \biggl[ \biggl( \int_{\mathbb{\{|\xi| \geq\mathrm{A} \}}}
|\Delta'|^2 \,\ud\xi\biggr)^{1/2} ; U^{c} \biggr]\nonumber
\end{eqnarray}
with
\[
\mathrm{A} = \mathrm{A}(\beta) := \frac{\sigma y_0}{ (
\sum_{j=1}^{\nu} \pi_{j}^{4} )^{1/4}} .
\]
At this stage, we apply (\ref{eq:restopunt4fast}) to the
evaluation of the first integral in (\ref{eq:somme}) after observing
that the function $\eta_{4,n}(\xi)$ here equals $e^{-\xi^2/2}$ almost
surely since $\kappa_4 = 0$. This leads to
%
%
\begin{eqnarray} \label{eq:primosomme}
&&\biggl( \int_{{ \{ |\xi| \leq\mathrm{A} \}}} |\Delta|^2 \,\ud
\xi\biggr)^{1/2} \nonumber\\
&&\qquad\leq2 \sqrt{2 \Gamma(17/2)} C_{4}^{\ast}
\Biggl( \sum_{j=1}^{\nu} \pi_{j}^6 \Biggr) \\
&&\qquad\quad{} + \sqrt{2} C_{4}^{\ast
} \Biggl[ \int_{\mathbb{R}} \xi^8 (1 + \xi
^4)^2 e^{-\xi^2} \Biggl( \sum_{j=1}^{\nu} \pi_{j}^4 \biggl| \tilde
{\epsilon}_4 \biggl(\frac{\pi_{j} \xi}{\sigma} \biggr) \biggr| \Biggr)^2
\,\ud\xi\Biggr]^{1/2}\nonumber
\end{eqnarray}
with
\[
\tilde{\epsilon}_4(x) :=
\cases{
\dfrac{\log\varphi_0(x) + (\sigma^2/2) x^2 - (\kappa_4/4!)
x^4}{x^4} , &\quad if $0 < |x| \leq\sigma y_0$, \cr
\tilde{\epsilon}_4(\sigma y_0) , &\quad if $ |x| > \sigma y_0$, \cr
0 , &\quad if $x = 0$.}
\]
Note that $\tilde{\epsilon}_4$ is a bounded continuous
function. Take expectations of both sides of (\ref{eq:primosomme}) and
recall (\ref{eq:formulaGR}) to obtain
%
%
\begin{eqnarray} \label{eq:primosommesperanza}\quad
&&\mathsf{E}_t \biggl( \int_{{ \{ |\xi| \leq\mathrm{A} \}}} |\Delta|^2
\,\ud
\xi\biggr)^{1/2}\nonumber\\
&&\qquad \leq2 \sqrt{2 \Gamma(17/2)} C_{4}^{\ast} e^{-{(3/8)}t} \\
&&\qquad\quad{} + \sqrt{2} C_{4}^{\ast} \mathsf{E}_t \Biggl[ \int_{\mathbb{R}} \xi
^8 (1 + \xi^4)^2 e^{-\xi^2} \Biggl( \sum_{j=1}^{\nu} \pi_{j}^4 \biggl|
\tilde{\epsilon}_4 \biggl(\frac{\pi_{j} \xi}{\sigma} \biggr) \biggr|
\Biggr)^2 \,\ud\xi\Biggr]^{1/2} .\nonumber
\end{eqnarray}
In view of Section \ref{app4},
%
%
\begin{equation} \label{eq:mattia}
\lim_{t \rightarrow+\infty} \rho_{0}^{(1)}(t) = 0 ,
\end{equation}
where
\[
\rho_{0}^{(1)}(t) := e^{{(1/4)}t} \mathsf{E}_t \Biggl[ \int
_{\mathbb{R}} \xi^8 (1 + \xi^4)^2 e^{-\xi^2} \Biggl( \sum_{j=1}^{\nu}
\pi
_{j}^4 \biggl| \tilde{\epsilon}_4 \biggl(\frac{\pi_{j} \xi}{\sigma}
\biggr) \biggr| \Biggr)^2 \,\ud\xi\Biggr]^{1/2} .
\]
Similarly, apply (\ref{eq:restopunt4der}) to evaluate the
second integral in (\ref{eq:somme}) as follows:
%
%
\begin{eqnarray}\label{eq:secondosomme}
&&\biggl( \int_{{ \{ |\xi| \leq\mathrm{A} \}}} |\Delta'|^2 \,\ud
\xi\biggr)^{1/2}\nonumber\\
&&\qquad \leq4 \sqrt{\Gamma(19/2)} C_{4}^{\ast} \Biggl(
\sum_{j=1}^{\nu} \pi_{j}^6 \Biggr) \nonumber\\[-8pt]\\[-8pt]
&&\qquad\quad{} + 2\sqrt{2} C_{4}^{\ast} \Biggl[ \int_{\mathbb{R}} \xi^6 (1 + \xi
^{12}) e^{-\xi^2} \Biggl( \sum_{j=1}^{\nu} \pi_{j}^4 \biggl| \tilde
{\epsilon}_4 \biggl(\frac{\pi_{j} \xi}{\sigma} \biggr) \biggr| \Biggr)^2
\,\ud\xi\Biggr]^{1/2} \nonumber\\
&&\qquad\quad{} + 2\sqrt{2} C_{4}^{\ast} \Biggl[ \int_{\mathbb{R}} \xi^6 (1 + \xi
^{12}) e^{-\xi^2} \Biggl( \sum_{j=1}^{\nu} \pi_{j}^4 \biggl| \tilde
{\varrho}_4 \biggl(\frac{\pi_{j} \xi}{\sigma} \biggr) \biggr| \Biggr)^2
\,\ud\xi\Biggr]^{1/2}\nonumber
\end{eqnarray}
with
\[
\tilde{\varrho}_4(x) :=
\cases{
x \dfrac{\ud}{\ud x}\tilde{\epsilon}_{4}(x) , &\quad if $0 < |x| <
\sigma y_0$, \vspace*{2pt}\cr
\displaystyle l := \lim_{u \uparrow\sigma y_0} \tilde{\varrho}_4(u), &\quad if
$|x| \geq\sigma y_0$, \vspace*{2pt}\cr
0 , &\quad if $x = 0$.}
\]
Once again, take expectations of both sides of (\ref
{eq:secondosomme}) and use (\ref{eq:formulaGR}) to get
%
%
\begin{eqnarray}\label{eq:secondosommesperanza}\qquad
&&\mathsf{E}_t \biggl( \int_{{ \{ |\xi| \leq\mathrm{A} \}}} |\Delta'|^2
\,\ud
\xi\biggr)^{1/2} \nonumber\\
&&\qquad\leq4 \sqrt{\Gamma(19/2)} C_{4}^{\ast} e^{-{(3/8)}t}
\nonumber\\[-8pt]\\[-8pt]
&&\qquad\quad{} + 2\sqrt{2} C_{4}^{\ast} \mathsf{E}_t \Biggl[ \int_{\mathbb{R}} \xi
^6 (1 + \xi^{12}) e^{-\xi^2} \Biggl( \sum_{j=1}^{\nu} \pi_{j}^4 \biggl|
\tilde{\epsilon}_4 \biggl(\frac{\pi_{j} \xi}{\sigma} \biggr) \biggr|
\Biggr)^2 \,\ud\xi\Biggr]^{1/2} \nonumber\\
&&\qquad\quad{} + 2\sqrt{2} C_{4}^{\ast} \mathsf{E}_t \Biggl[ \int_{\mathbb{R}} \xi
^6 (1 + \xi^{12}) e^{-\xi^2} \Biggl( \sum_{j=1}^{\nu} \pi_{j}^4 \biggl|
\tilde{\varrho}_4 \biggl(\frac{\pi_{j} \xi}{\sigma} \biggr) \biggr|
\Biggr)^2 \,\ud\xi\Biggr]^{1/2} .\nonumber
\end{eqnarray}
Another application of Section \ref{app4} leads us to state the following
important facts:
%
%
\begin{equation} \label{eq:paolo}
\lim_{t \rightarrow+\infty} \rho_{0}^{(2)}(t) = \lim_{t
\rightarrow
+\infty} \rho_{0}^{(3)}(t) = 0 ,
\end{equation}
where
\[
\rho_{0}^{(2)}(t) := e^{{(1/4)}t} \mathsf{E}_t \Biggl[ \int
_{\mathbb{R}} \xi^6 (1 + \xi^{12}) e^{-\xi^2} \Biggl( \sum_{j=1}^{\nu}
\pi_{j}^4 \biggl| \tilde{\epsilon}_4 \biggl(\frac{\pi_{j} \xi}{\sigma}
\biggr) \biggr| \Biggr)^2 \,\ud\xi\Biggr]^{1/2}
\]
and
\[
\rho_{0}^{(3)}(t) := e^{{(1/4)}t} \mathsf{E}_t \Biggl[ \int
_{\mathbb{R}} \xi^6 (1 + \xi^{12}) e^{-\xi^2} \Biggl( \sum_{j=1}^{\nu}
\pi_{j}^4 \biggl| \tilde{\varrho}_4 \biggl(\frac{\pi_{j} \xi}{\sigma}
\biggr) \biggr| \Biggr)^2 \,\ud\xi\Biggr]^{1/2} .
\]
After determining upper bounds for integrals of the type
$\int_{\{|\xi| \leq\mathrm{A}\}}$, it remains to examine the remaining
summands in (\ref{eq:somme}). Minkowski's inequality yields
\[
\biggl( \int_{{\mathbb{\{|\xi| \geq\mathrm{A} \}}}} |\Delta|^2 \,\ud
\xi\biggr)^{1/2} \leq\biggl( \int_{{\mathbb{\{|\xi| \geq
\mathrm{A} \}}}} |\varphi^{\ast}(\xi/\sigma)|^2 \,\ud\xi\biggr)^{1/2}
+ \biggl( \int_{{\mathbb{\{|\xi| \geq\mathrm{A} \}}}} |e^{- {
\xi^2}/{2}}|^2 \,\ud\xi\biggr)^{1/2}
\]
and
\begin{eqnarray*}
&&\biggl( \int_{{\mathbb{\{|\xi| \geq\mathrm{A} \}}}} |\Delta'|^2 \,\ud
\xi\biggr)^{1/2} \\
&&\qquad\leq\biggl( \int_{{\mathbb{\{|\xi| \geq
\mathrm{A} \}}}} \biggl|\frac{\ud}{\ud\xi}\varphi^{\ast}(\xi/\sigma
)\biggr|^2 \,\ud
\xi\biggr)^{1/2} + \biggl( \int_{{\mathbb{\{|\xi| \geq\mathrm{A}
\}}}} |\xi e^{- {\xi^2}/{2}}|^2 \,\ud\xi\biggr)^{1/2} .
\end{eqnarray*}
From a well-known inequality, proved in, for example, Lemma 2
of VII.1 in Feller (\citeyear{F68}), and since $\max_{x \geq0} x^k e^{- \alpha
x^2} = [k/(2 e \alpha)]^{k/2}$, one obtains
\[
\biggl( \int_{{\mathbb{\{|\xi| \geq\mathrm{A} \}}}} e^{- \xi^2} \,\ud
\xi\biggr)^{1/2} \leq\biggl( \frac{15}{2} \biggr)^{15/4} e^{-15/4}
(\sigma y_0)^{-8} \sum_{j=1}^{\nu} \pi_{j}^{6}
\]
and
\[
\biggl( \int_{{\mathbb{\{|\xi| \geq\mathrm{A} \}}}} \xi^2 e^{- \xi^2}
\,\ud\xi\biggr)^{1/2} \leq\frac{2 + \sqrt{2}}{2} \biggl( \frac
{17}{2} \biggr)^{17/4} e^{-15/4} (\sigma y_0)^{-8} \sum_{j=1}^{\nu} \pi
_{j}^{6} .
\]
Equation (\ref{eq:formulaGR}) can then be applied to obtain
%
%
\begin{equation} \label{eq:della}
\mathsf{E}_t \biggl( \int_{{\mathbb{\{|\xi| \geq\mathrm{A} \}}}} e^{-
\xi
^2} \,\ud
\xi\biggr)^{1/2} \leq\biggl( \frac{15}{2} \biggr)^{15/4} e^{-15/4}
(\sigma y_0)^{-8} e^{-{(3/8)}t}
\end{equation}
and
%
%
\begin{equation} \label{eq:croce}\quad
\mathsf{E}_t \biggl( \int_{{\mathbb{\{|\xi| \geq\mathrm{A} \}}}} \xi^2
e^{- \xi^2}
\,\ud\xi\biggr)^{1/2} \leq\frac{2 + \sqrt{2}}{2} \biggl( \frac
{17}{2} \biggr)^{17/4} e^{-15/4} (\sigma y_0)^{-8} e^{-{(3/8)}t} .
\end{equation}
At this point, to control the remaining integrals over $\{
|\xi| \geq\mathrm{A} \}$, we proceed as in formula (30) of DGR to write
%
%
\begin{eqnarray} \label{eq:ibra}\quad
&&\biggl[ \biggl( \int_{{\mathbb{\{|\xi| \geq\mathrm{A} \}}}}
|\varphi^{\ast}(\xi/\sigma)|^2 \,\ud\xi\biggr)^{1/2} + \biggl(
\int_{{\mathbb{\{|\xi| \geq\mathrm{A} \}}}} \biggl|\frac{\ud}{\ud\xi}
\varphi^{\ast}(\xi/\sigma)\biggr|^2 \,\ud\xi\biggr)^{1/2}
\biggr] \cdot\mathbh{1}_{U^{c}} \nonumber\\[-8pt]\\[-8pt]
&&\qquad\leq2\sqrt{2} \biggl( \int
_{\mathrm
{A}}^{+\infty}
|\varphi^{\ast}(\xi/\sigma)| \,\ud\xi\biggr)^{1/2} \cdot
\mathbh{1}_{U^{c}} + \sqrt{2|\varphi^{\ast}(\mathrm{A}/\sigma)|}
\cdot\mathbh{1}_{U^{c}}.\nonumber
\end{eqnarray}
For $\overline{\omega}$ in $U^c$, the bound
\[
\mathrm{A}(\overline{\omega}) \leq\frac{\sigma^3}{2 \overline
{\mathfrak{m}}_3 \sum_{j=1}^{\nu(\overline{\omega})}
|\pi_{j}(\overline{\omega})|^3}
\]
holds true, thanks to the definition of $\overline{\delta}$
and the Lyapunov inequality. Thus, Lemma~12 in Chapter 6 of Petrov
(\citeyear{P75}) can be applied to the characteristic function $\varphi^{\ast
}(\xi
/\sigma)$ with $b = 1/2$ to deduce
\begin{eqnarray*}
\sqrt{2 |\varphi^{\ast}(\mathrm{A}/\sigma)| } &\leq& \sqrt{2} e^{-
A^2/12 } \leq\sqrt{2} (48/e)^{4} \mathrm
{A}^{-8} \\
&=& \sqrt{2} (48/e)^{4} (\sigma y_0)^{-8}
\sum_{j=1}^{\nu} \pi_{j}^{6},
\end{eqnarray*}
which entails
that
%
%
\begin{equation} \label{eq:esternoPetrov}
\mathsf{E}_t \sqrt{2 | \varphi^{\ast}(\mathrm{A}/\sigma) | } \leq
\sqrt{2} (48/e)^{4} (\sigma y_0)^{-8}
e^{-{(3/8)}t} .
\end{equation}
It remains to analyze
\[
\biggl( \int_{\mathrm{A}}^{+\infty} | \varphi^{\ast}(\xi/\sigma)| \,
\ud
\xi
\biggr)^{1/2} \cdot\mathbh{1}_{U^{c}} = \Biggl( \int_{\mathrm{A}}^{+\infty}
\prod_{j=1}^{\nu} \biggl| \varphi_0\biggl(\frac{\pi_j \xi}{\sigma}\biggr)
\biggr| \,\ud\xi\Biggr)^{1/2} \cdot\mathbh{1}_{U^{c}} .
\]
An estimate of this term is made using Proposition 2.2 in
DGR, together with (33), (34) and (35) therein, with $\overline
{\varepsilon} = 1/(2\overline{n}!)$. We then have
%
%
\begin{eqnarray}\label{eq:code}
\biggl( \int_{\mathrm{A}}^{+\infty} | \varphi^{\ast}(\xi/\sigma)| \,
\ud
\xi
\biggr)^{1/2} \cdot\mathbh{1}_{U^{c}} &\leq& \biggl[ \lambda\sigma\int_{\mathrm
{A}/\lambda\sigma}^{+\infty}
\biggl(\frac{1}{\overline{\varepsilon} \eta^{2\overline{n}}}
\biggr)^{\alpha} \,\ud\eta\biggr]^{1/2} \nonumber\\[-8pt]\\[-8pt]
&=& D
\Biggl( \sum_{j=1}^{\nu} \pi_{j}^{4} \Biggr)^{ ({ 2 \overline{n}
\alpha- 1})/{8} } .\nonumber
\end{eqnarray}
The definition of $\overline{n}$ in (\ref{eq:insiemecattivo})
yields $(2 \alpha\overline{n} -
1)/8 = 2$. Moreover,
%
%
\begin{eqnarray} \label{eq:ultimacostante}
D :\!&=& \frac{1}{4 \overline{\varepsilon}^{\alpha/2}} \frac{(\lambda
\sigma)^{17/2}}{(\sigma y_0)^8} \nonumber\\[-8pt]\\[-8pt]
&\leq&2^{13/4} \biggl[ \biggl(
\frac{3}{2\sigma^2} \biggr)^{17/4} +\biggl ( \frac{2}{1 -
M} \biggr)^{17/4} (L_p)^{17/2p} \biggr]\nonumber
\end{eqnarray}
with
\[
L_p := \sup_{\xi\in\mathbb{R}} [ |\xi|^p \cdot|\varphi_0(\xi)| ]
\]
and
\[
M = \exp\biggl\{ - \frac{3\pi^2}{64
( 3 + (L_p)^{4/p} )^2} \biggl(\frac{\sqrt{2}\sigma}{8
\lceil2/p \rceil\sigma^3 + 40\pi\sqrt{\lceil2/p \rceil\mathfrak{m}_4}}
\biggr)^2 \biggr\} .
\]
Taking expectation in (\ref{eq:code}) gives
%
%
\begin{equation} \label{eq:esternocaratteristica}
\mathsf{E}_t \biggl[ \biggl( \int_{\mathrm{A}}^{+\infty} |
\varphi^{\ast}(\xi/\sigma)| \,\ud\xi\biggr)^{1/2} \cdot\mathbh{1}_{U^{C}}
\biggr] \leq D e^{-{(3/8)}t} .
\end{equation}
The claimed upper bound (\ref{eq:D+delta}) follows from (\ref
{eq:stimasulbrutto}), (\ref{eq:primosommesperanza}), (\ref
{eq:secondosommesperanza}), (\ref{eq:della}), (\ref{eq:croce}), (\ref
{eq:esternoPetrov}) and~(\ref{eq:esternocaratteristica}).

\subsection{Proof of Theorems \protect\ref{thm:D+delta} and \protect
\ref{thm:piuottimo} when $2 \chi+ \delta> 4$}\label{sec44}

This proof differs from the previous one only in the choice of the
constants. One can start from (\ref{eq:primastima}) under hypothesis
(\ref{hypoH}). Thanks to (\ref{hypoH}) and the hypotheses of the theorems to be proven, one
can apply Proposition 2.2 of DGR to get (\ref{eq:magggamma}) with
$\alpha= (2 \cdot\lceil2/p \rceil)^{-1}$.

Now, define $U$ exactly as in (\ref{eq:insiemecattivo}) with
$\overline{n} = [k(k + 2) + 1] \cdot\lceil2/p \rceil$ and
\[
\overline{\delta} = \min\biggl\{ \frac{1}{2^{\overline{n}} \overline
{n}!}; \frac{\sigma^8}{16 y_{0}^{4} \overline{\mathfrak{m}}_{3}^{4}}
\biggr\} \leq\frac{1}{2^{\overline{n}} \overline{n}!} .
\]
The probability of $U$ is then estimated, under each $\mathsf
{P}_t$, using the facts that
\[
\mathsf{P}_t \{ \nu\leq\overline{n} \} \leq\overline{n} e^{-t}
\quad\mbox{and}\quad \mathsf{P}_t \Biggl\{ \prod_{j=1}^{\nu
} \pi_j = 0 \Biggr\} = 0 ,
\]
whereas, for the third component of the union in the
definition of $U$, one can combine Markov's (with exponent $k/2$) and
Lyapounov's inequalities to get
\begin{eqnarray*}
\mathsf{P}_t \Biggl\{ \sum_{j=1}^{\nu} \pi_{j}^{4} \geq\overline{\delta}
\Biggr\} &\leq& \frac{1}{\overline{\delta}{}^{k/2}} \mathsf{E}_t \Biggl[
\Biggl( \sum_{j=1}^{\nu} \pi_{j}^{4} \Biggr)^{k/2} \Biggr] \\
&\leq&\frac{1}{\overline{\delta}{}^{k/2}} \mathsf{E}_t \Biggl[ \sum
_{j=1}^{\nu} \pi_{j}^{k+2} \Biggr] \leq\frac{1}{\overline{\delta
}^{k/2}} e^{- (1 - 2\alpha_{k+2}) t} .
\end{eqnarray*}
Thus,
%
%
\begin{equation} \label{eq:probinsiemecattivok}
\mathsf{P}_t (U) \leq[\overline{n} + 1/\overline{\delta}{}^{k/2}]
e^{- (1 - 2\alpha_{k+2}) t}\qquad (t \geq0) .
\end{equation}
Now, split the term $\mathsf{E}_t [ \Vert\frac{\ud
}{\ud v} \mathsf{F}^{\ast}(\sigma v) - g_{1}(v) \Vert_{1}]$ into
the sum of two contributions, exactly as in (\ref{eq:secondastima}),
and note that (\ref{eq:probinsiemecattivok}) entails
that
%
%
\begin{equation} \label{eq:stimasulbruttok}\qquad
\mathsf{E}_t \biggl[ \biggl\Vert\frac{\ud}{\ud v}
\mathsf{F}^{\ast}(\sigma v) - g_{1}(v) \biggr\Vert_{1}; U \biggr] \leq2\mathsf
{P}_t (U) \leq2 [\overline{n} + 1/\overline{\delta
}^{k/2}] e^{- (1 - 2\alpha_{k+2}) t} .
\end{equation}
To control the integral over $U^c$, we once again invoke
Beurling's inequality (see Proposition 4.1 in DGR) to write (\ref
{eq:beurlingsulbrutto}). Applicability of this result rests on the same
arguments as those provided in Section \ref{sec43}. The right-hand side of
(\ref
{eq:beurlingsulbrutto}) is split into a sum of four terms, exactly as
in (\ref{eq:somme}), with
\[
\mathrm{A} = \mathrm{A}(\beta) := \frac{\sigma y_0}{ (
\sum_{j=1}^{\nu} \pi_{j}^{4} )^{1/(k + \delta)}} .
\]
Now, apply (\ref{eq:restoint}) to the evaluation of the first
integral in (\ref{eq:somme}), noting that the function $\eta
_{k,n}(\xi
)$ equals $e^{-\xi^2/2}$ almost surely since $\kappa_{2r} = 0$ for $r =
2, \ldots, \chi$. This leads to
%
%
\begin{equation} \label{eq:internok}
\mathsf{E}_t \biggl[ \biggl( \int_{{\mathbb{\{|\xi| \leq\mathrm{A}
\}}}} |\Delta|^2 \,\ud\xi\biggr)^{1/2} \biggr] \leq C_{k,\delta
}^{\ast} a_k \cdot e^{- (1 - 2\alpha_{2\chi+ \delta}) t}
\end{equation}
and
%
%
\begin{equation} \label{eq:derivatainternok}
\mathsf{E}_t \biggl[ \biggl( \int_{{\mathbb{\{|\xi| \leq\mathrm{A}
\}}}} |\Delta'|^2 \,\ud\xi\biggr)^{1/2} \biggr] \leq C_{k,\delta
}^{\ast} a_k \cdot e^{- (1 - 2\alpha_{2\chi+ \delta}) t}.
\end{equation}
After determining upper bounds for integrals of the type $\int
_{\{|\xi| \leq\mathrm{A}\}}$, it remains to examine the remaining
summands in (\ref{eq:somme}). Minkowski's inequality gives
\[
\biggl( \int_{{\mathbb{\{|\xi| \geq\mathrm{A} \}}}} |\Delta|^2 \,\ud
\xi\biggr)^{1/2} \leq\biggl( \int_{{\mathbb{\{|\xi| \geq
\mathrm{A} \}}}} |\varphi^{\ast}(\xi/\sigma)|^2 \,\ud\xi\biggr)^{1/2}
+ \biggl( \int_{{\mathbb{\{|\xi| \geq\mathrm{A} \}}}} |e^{- {
\xi^2}/{2}}|^2 \,\ud\xi\biggr)^{1/2}
\]
and
\begin{eqnarray*}
\biggl( \int_{{\mathbb{\{|\xi| \geq\mathrm{A} \}}}} |\Delta'|^2 \,\ud
\xi\biggr)^{1/2}
&\leq&\biggl( \int_{{\mathbb{\{|\xi| \geq
\mathrm{A} \}}}} \biggl|\frac{\ud}{\ud\xi}\varphi^{\ast}(\xi/\sigma
)\biggr|^2 \,\ud
\xi\biggr)^{1/2}\\
&&{} + \biggl( \int_{{\mathbb{\{|\xi| \geq\mathrm{A}
\}}}} |\xi e^{- {\xi^2}/{2}}|^2 \,\ud\xi\biggr)^{1/2} .
\end{eqnarray*}
Integrals involving the Gaussian density are controlled as in
the previous subsection, giving
%
%
\begin{eqnarray} \label{eq:esternogaussk}
&&\mathsf{E}_t \biggl( \int_{{\mathbb{\{|\xi| \geq\mathrm{A} \}}}} e^{-
\xi^2} \,\ud\xi\biggr)^{1/2}\nonumber\\[-8pt]\\[-8pt]
&&\qquad \leq\biggl(\frac{k (k + 2) - 1}{2 e}
\biggr)^{({k (k + 2) - 1})/{4}}
(\sigma y_0)^{- k (k + 2)/2} e^{- (1 - 2\alpha_{k+2}) t}\nonumber
\end{eqnarray}
and
%
%
\begin{eqnarray}\label{eq:esternoderivatagaussk}\quad
&&\mathsf{E}_t \biggl( \int_{{\mathbb{\{|\xi| \geq\mathrm{A} \}}}}
\xi^2 e^{- \xi^2} \,\ud\xi\biggr)^{1/2} \nonumber\\[-8pt]\\[-8pt]
&&\qquad\leq\frac{2 + \sqrt{2}}{2} \biggl(\frac{k (k + 2) + 1}{2 e}
\biggr)^{({k (k + 2) + 1})/{4}}
(\sigma y_0)^{- k (k + 2)/2} e^{- (1 - 2\alpha_{k+2}) t} .\nonumber
\end{eqnarray}
To control the remaining integrals over the region $\{ |\xi|
\geq\mathrm{A} \}$, we proceed as before, writing (\ref
{eq:ibra}). For $\overline{\omega}$ in $U^c$, the bound
\[
\mathrm{A}(\overline{\omega}) \leq\frac{\sigma^3}{2 {\overline
{\mathfrak{m}}_3 \sum_{j=1}^{\nu(\overline{\omega})}}
|\pi_{j}(\overline{\omega})|^3}
\]
holds true, thanks to the definition of $\overline{\delta}$
and the Lyapunov inequality. We then set $b = 1/2$ in Lemma 12 from
Chapter 6 of Petrov (\citeyear{P75}) to deduce
that
\begin{eqnarray*}
&&\sqrt{2 |\varphi^{\ast}(\mathrm{A}/\sigma)| } \\
&&\qquad\leq\sqrt{2} e^{-
{{A}^2/12} } \\
&&\qquad\leq\sqrt{2} \biggl( \frac{3 k (k + 2)}{e} \biggr)^{({k (k +2)})/{4}}
(\sigma y_0)^{- ({k (k +2)})/{2}}
\cdot\Biggl( \sum_{j=1}^{\nu} \pi_{j}^4 \Biggr)^{({k (k + 2)})/({2 (k
+ \delta)})} \\
&&\qquad\leq\sqrt{2} \biggl( \frac{3 k (k + 2)}{e} \biggr)^{({k (k +2)})/{4}}
(\sigma y_0)^{- ({k (k +2)})/{2}} \cdot\Biggl( \sum_{j=1}^{\nu} \pi
_{j}^{k + 2} \Biggr)
\end{eqnarray*}
and, therefore,
%
%
\begin{eqnarray}\label{eq:esternoPetrovk}
\mathsf{E}_t \sqrt{2 | \varphi^{\ast}(\mathrm{A}/\sigma) | }
&\leq&
\sqrt{2} \biggl( \frac{3 k (k + 2)}{e} \biggr)^{({k (k +2)})/{4}}
(\sigma y_0)^{- ({k (k +2)})/{2}} \nonumber\\[-8pt]\\[-8pt]
&&{}\times e^{-(1 - 2\alpha_{k+2}) t} .\nonumber
\end{eqnarray}

Finally, in regard to $ ( \int_{\mathrm{A}}^{+\infty} | \varphi
^{\ast
}(\xi/\sigma)| \,\ud\xi
)^{1/2} \cdot\mathbh{1}_{U^{c}}$, one can write
%
%
\begin{equation} \label{eq:codek}
\biggl( \int_{\mathrm{A}}^{+\infty} | \varphi^{\ast}(\xi/\sigma)| \,
\ud
\xi
\biggr)^{1/2} \cdot\mathbh{1}_{U^{c}} = D_k
\Biggl( \sum_{j=1}^{\nu} \pi_{j}^{4} \Biggr)^{ ({ 2 \overline{n}
\alpha- 1})/{8} },
\end{equation}
where the constant $D_k$ is given by
\[
\sqrt{\frac{\lambda\sigma(2\overline{n}!)^{\alpha}}{2\alpha
\overline
{n} - 1}} \biggl(\frac{\lambda}{y_0} \biggr)^{({2\alpha\overline{n}
- 1})/{2}} .
\]

The definition of $\overline{n}$ given at the beginning of this
subsection yields $(2 \alpha\overline{n} - 1)/8 > k/2$. Now, taking
expectation in (\ref{eq:codek}) entails
that
%
%
\begin{equation} \label{eq:esternocaratteristicak}
\mathsf{E}_t \biggl[ \biggl( \int_{A}^{+\infty} |
\varphi^{\ast}(\xi/\sigma)| \,\ud\xi\biggr)^{1/2} \cdot\mathbh{1}_{U^{C}}
\biggr] \leq D_k e^{- (1 - 2\alpha_{k+2})t} .
\end{equation}
To obtain (\ref{eq:D+k+delta}), it will suffice to combine
the previous inequalities.

\begin{appendix}\label{app}
\section*{Appendix}

This appendix contains all of the elements which are necessary to
complete the proofs given in Section \ref{sec4}. It is split into four parts.
The first focuses on a quantification of the numbers $y_0$ such that
the Fourier--Stieltjes transform of a symmetric probability law turns
out to be greater than $1/2$ on $[-y_0, y_0]$. The second presents the
proofs of Lemmas \ref{lm:approssimazioni4} and \ref
{lm:approssimazioni}. The third aims to clarify the conclusion of the
proof of Proposition \ref{prop:menouno}. Finally, the fourth provides a
proof for (\ref{eq:mattia}) and (\ref{eq:paolo}).

\subsection{Specification of $y_0$}\label{app1}
\textit{Let $\psi$ be the
Fourier--Stieltjes transform of a symmetric probability law
$\zeta$ on $(\mathbb{R}, \mathscr{B}(\mathbb{R}))$},
\textit{namely $\psi(\xi) := \int_{\mathbb{R}} e^{i \xi x} \zeta(\ud x)$
for every real $\xi$. Assume that $\mathfrak{m}_4 := \int
_{\mathbb{R}} x^4 \zeta(\ud x)$ is finite\vspace*{2pt} and put
$\sigma
^2 := \int_{\mathbb{R}} x^2 \zeta(\ud x)$, $y_0 := \{ [-6 \sigma^2 +
(36 \sigma^4 + 12\mathfrak{m}_4)^{1/2}]/\mathfrak{m}_4 \}^{1/2}$.
If $|\xi| \leq y_0$}, \textit{then $\psi(\xi) \geq1/2$.}
\begin{pf}
By the Taylor expansion for characteristic
functions, one can write $\psi(\xi) = 1 - (\sigma^2/2) \xi^2 +
R(\xi)$
with $|R(\xi)| \leq(\mathfrak{m}_4/24) \xi^4$; see, for example,
Section~8.4 in Chow and Teicher (\citeyear{CT97}). The desired bound is obtained if
\[
1 - \frac{\sigma^2}{2} \xi^2 - \frac{\mathfrak{m}_4}{24} \xi^4
\geq\frac{1}{2}
\]
holds true for every $\xi$ belonging to some interval. Now,
one can note that the biquadratic equation $\mathfrak{m}_4 \xi^4 + 12
\sigma^2 \xi^2 - 12 = 0$ possesses exactly two real solutions, namely
$\pm y_0$, and the previous inequality is satisfied for every $\xi$ in
$[-y_0, y_0]$.
\end{pf}

\subsection{Proofs of Lemmas \protect\ref{lm:approssimazioni4} and \protect\ref
{lm:approssimazioni}}\label{app2}

\mbox{}

\begin{pf*}{Proof of Lemma \protect\ref{lm:approssimazioni4}}
Set $\psi_{j,n}$ for the characteristic
function of $Y_{j,n}$ $(j = 1, 2, \ldots, n)$ and use the definition of
$V_n$, combined with independence, to write
\[
\psi_n(\xi) = \prod_{j=1}^n \psi_{j,n}(\xi) = \prod_{j=1}^n \psi\biggl(
\frac{c_{j,n} \xi}{\sigma} \biggr) .
\]
If $|\xi| \leq A_{4,n}$, then it easily follows that
\[
\biggl| \frac{c_{j,n} \xi}{\sigma} \biggr| \leq\Biggl| \frac
{c_{j,n} \sigma y_0}{\sigma} \Biggl( \sum_{r=1}^n c_{r,n}^4
\Biggr)^{-1/4} \Biggr| \leq y_0 .
\]
Now, using elementary properties of the logarithm, one can
combine expansion (\ref{eq:logtaylor}) with property (\ref
{eq:sferadiscreta}) of each array $\{ c_{1,n}, \ldots, c_{n,n} \}$ to obtain
\begin{eqnarray*}
\log\psi_n(\xi) &=& \sum_{j=1}^n \log\psi_{j,n}(\xi) \\
&=& \sum_{j=1}^n \biggl[ -\frac{1}{2} \sigma^2 \frac{c_{j,n}^2 \xi
^2}{\sigma^2} + \frac{1}{4!} \kappa_4 \frac{c_{j,n}^4 \xi
^4}{\sigma^4}
+ \frac{c_{j,n}^4 \xi^4}{\sigma^4} \epsilon_4 \biggl(\frac{c_{j,n} \xi
}{\sigma} \biggr) \biggr] \\
&=& -\frac{1}{2} \xi^2 + \frac{\tilde{\lambda}_{2,n}}{4!} \xi^4 +
R_4(\xi),
\end{eqnarray*}
where
\[
R_4(\xi) := \sum_{j=1}^n \frac{c_{j,n}^4 \xi^4}{\sigma^4} \epsilon_4
\biggl(\frac{c_{j,n} \xi}{\sigma} \biggr) .
\]
Inverting the logarithm, one gets
%
%
\begin{equation} \label{eq:petrov1}
\psi_n(\xi) = e^{- \xi^2/2} \cdot\exp\biggl\{ \frac{\tilde{\lambda
}_{2,n}}{4!} \xi^4 \biggr\} \cdot\exp\{ R_4(\xi) \} .
\end{equation}
It is easily verified that the restrictions $|u| := |\tilde
{\lambda}_{2,n} \xi^4|/4! \leq\kappa_4 y_0/4!$ and $|R_4(\xi)|
\leq
M_{0}^{(4)} y_{0}^{4}$ hold true when $|\xi| \leq A_{4,n}$, and that
$\tilde{\lambda}_{2,n} \xi^4/4! = \tilde{P}_{1,n}(\xi)$. Finally, set
$F(x) := e^x - 1 - x$. At this point, we have all the tools needed to
prove (\ref{eq:restopunt4}) and (\ref{eq:restopunt4fast}). Indeed,
\begin{eqnarray*}
|\psi_n(\xi) - \eta_{4,n}(\xi)| &=& e^{- \xi^2/2} | e^u \exp\{
R_4(\xi) \} - 1 - u | \\
&=& e^{- \xi^2/2} | e^u \exp\{ R_4(\xi) \} - e^u + F(u) |
\\
&\leq& e^{- \xi^2/2} e^u | \exp\{ R_4(\xi) \} - 1 | + e^{- \xi^2/2}
|F(u)| .
\end{eqnarray*}
By elementary arguments, if $x$ is any real number satisfying
$|x| \leq c$, one has
\[
| e^x - 1| \leq e^{|x|} - 1 \leq\biggl( \frac{e^c -1}{c}
\biggr) |x| .
\]
This fact can be applied to $R_4(\xi)$ to get
\[
| \exp\{ R_4(\xi) \} - 1 | \leq\xi^4 \cdot\biggl( \frac
{e^{M_{0}^{(4)} y_{0}^{4}} -1}{\sigma^4 y_{0}^{4}} \biggr) \cdot\Biggl(
\sum_{j=1}^n c_{j,n}^4 \biggl| \epsilon_4 \biggl(\frac{c_{j,n} \xi
}{\sigma} \biggr) \biggr| \Biggr) .
\]
Moreover, since the inequality
\[
|F(u)| \leq\max_{|x| \leq\kappa_4 y_0/4!} \biggl[\biggl| \frac
{F(x)}{x^2}\biggr| \biggr] \xi^8 \Biggl( \sum_{j=1}^n c_{j,n}^4
\Biggr)^2
\]
holds, one can conclude that
%
%
\begin{eqnarray} \label{eq:stoppa}
&&|\psi_n(\xi) - \eta_{4,n}(\xi)| \nonumber\\
&&\qquad\leq e^{- \xi^2/2} \xi^4 \cdot\exp\biggl\{ \frac{\kappa_4 y_{0}^{4}}{4!}
\biggr\} \biggl( \frac{e^{M_{0}^{(4)} y_{0}^{4}} -1}{\sigma^4 y_{0}^{4}}
\biggr) \cdot\Biggl( \sum_{j=1}^n c_{j,n}^4 \biggl| \epsilon_4
\biggl(\frac{c_{j,n} \xi}{\sigma} \biggr) \biggr| \Biggr) \\
&&\qquad\quad{} + e^{- \xi^2/2} \max_{|x| \leq\kappa_4 y_{0}^{4}/4!} \biggl[\biggl|
\frac{F(x)}{x^2} \biggr| \biggr] \xi^8 \Biggl( \sum_{j=1}^n c_{j,n}^4
\Biggr)^2.\nonumber
\end{eqnarray}
After setting
\[
C^{\ast\ast}_4 := \exp\biggl\{ \frac{\kappa_4 y_{0}^{4}}{4!} \biggr\} \biggl( \frac
{e^{M_{0}^{(4)} y_{0}^{4}} -1}{\sigma^4 y_{0}^{4}} \biggr) + \max_{|x|
\leq\kappa_4 y_{0}^{4}/4!} \biggl[\biggl| \frac{F(x)}{x^2} \biggr| \biggr],
\]
the derivation of (\ref{eq:restopunt4}) and (\ref
{eq:restopunt4fast}) follows by rewriting (\ref{eq:stoppa})
in a more convenient form. To get (\ref{eq:restopunt4}), it is enough
to observe that $\sum_{j=1}^n c_{j,n}^4 \leq1$, while to deduce (\ref
{eq:restopunt4fast}), one can combine the inequality $ (\sum
_{j=1}^n c_{j,n}^4 )^2 \leq\sum_{j=1}^n c_{j,n}^6$ with $\max\{1
; \xi^4\} \leq(1 + \xi^4)$.

To prove (\ref{eq:restopunt4der}), we start from (\ref{eq:petrov1}) and
take the derivative with respect to $\xi$. Thus, one obtains
\begin{eqnarray*}
&&|\psi_{n}^{\prime} (\xi) - \eta_{4,n}^{\prime}(\xi)| \\
&&\qquad\leq\exp\{R_4(\xi)\} \cdot|R_{4}^{\prime}(\xi)| \cdot|\eta
_{4,n}(\xi
) +
F(u)e^{- \xi^2 /2}|\\
&&\qquad\quad{} + |\eta_{4,n}^{\prime}(\xi)| \cdot| \exp\{
R_4(\xi) \} - 1 | \\
&&\qquad\quad{} + \exp\{R_4(\xi)\} \cdot\biggl| \frac{\ud}{\ud\xi} F(u) \biggr|
\cdot e^{- \xi^2 /2} + \exp\{R_4(\xi)\} \cdot|F(u)| \cdot|\xi|
e^{- \xi^2 /2} .
\end{eqnarray*}
Arguing as in the first part of this proof,
we have
%
%
\begin{eqnarray}\label{eq:pippo1}
&&|\eta_{4,n}^{\prime}(\xi)| \cdot| \exp\{ R_4(\xi) \} - 1 |
\nonumber\\
&&\qquad\leq\biggl( \frac{e^{M_{0}^{(4)} y_{0}^{4}} -1}{\sigma^4 y_{0}^{4}}
\biggr) \cdot\biggl( 1 + \frac{\kappa_4}{4! \sigma^4} \biggr) |\xi|^5
(1 + \xi^4) e^{- \xi^2 /2} \\
&&\qquad\quad{}\times\Biggl( \sum_{j=1}^n c_{j,n}^4 \biggl|
\epsilon_4 \biggl(\frac{c_{j,n} \xi}{\sigma} \biggr) \biggr| \Biggr)\nonumber
\end{eqnarray}
and
%
%
\begin{eqnarray} \label{eq:pippo2}
&&\exp\{R_4(\xi)\} \cdot|F(u)| \cdot|\xi| e^{- \xi^2 /2}
\nonumber\\[-8pt]\\[-8pt]
&&\qquad\leq\max_{|x| \leq\kappa_4 y_0/4!} \biggl[\biggl| \frac{F(x)}{x^2}
\biggr| \biggr] \exp\bigl\{ M_{0}^{(4)} y_{0}^{4}\bigr\} |\xi|^9 e^{- \xi^2 /2}
\Biggl( \sum_{j=1}^n c_{j,n}^4 \Biggr)^2.\nonumber
\end{eqnarray}
Moreover,
%
%
\begin{eqnarray} \label{eq:pippo3}
&&\exp\{R_4(\xi)\} \cdot|R_{4}^{\prime}(\xi)| \cdot|\eta_{4,n}(\xi
) +
F(u)e^{- \xi^2 /2}| \nonumber\\
&&\qquad= \exp\{R_4(\xi)\} \cdot|R_{4}^{\prime}(\xi)|
\cdot
e^{- \xi^2 /2} e^u \nonumber\\[-8pt]\\[-8pt]
&&\qquad\leq\exp\bigl\{ M_{0}^{(4)} y_{0}^{4}\bigr\} \cdot\exp\biggl\{ \frac{\kappa_4
y_{0}^{4}}{4!}\biggr\} 4 \sigma^{-4} |\xi|^3 e^{- \xi^2 /2} \cdot
\nonumber\\
&&\qquad\quad{}\times\Biggl[ \sum_{j=1}^n c_{j,n}^4 \biggl| \epsilon_4 \biggl(\frac
{c_{j,n} \xi}{\sigma} \biggr) \biggr| + \sum_{j=1}^n c_{j,n}^4 \biggl|
\varrho_4 \biggl(\frac{c_{j,n} \xi}{\sigma} \biggr) \biggr|
\Biggr]\nonumber
\end{eqnarray}
and
\[
\frac{\ud}{\ud\xi} F(u) = \frac{\tilde{\lambda}_{2,n}}{3!} \xi
^3 (e^u - 1) ,
\]
whence
%
%
\begin{eqnarray} \label{eq:pippo4}\quad
&&\exp\{R_4(\xi)\} \cdot\biggl| \frac{\ud}{\ud\xi} F(u) \biggr| \cdot
e^{- \xi^2 /2} \nonumber\\[-8pt]\\[-8pt]
&&\qquad\leq\exp\bigl\{ M_{0}^{(4)} y_{0}^{4}\bigr\} \frac{\kappa_{4}^{2}}{3! 4!
\sigma^8} \biggl( \frac{\exp\{ {\kappa_4 y_{0}^{4}}/{4!} \} -
1}{{\kappa_4 y_{0}^{4}}/{4!}} \biggr) |\xi|^7 e^{- \xi^2 /2} \Biggl( \sum
_{j=1}^n c_{j,n}^4 \Biggr)^2 .\nonumber
\end{eqnarray}
Now, set
\begin{eqnarray*}
C_{4}^{\ast\ast\ast} &:=& \biggl( \frac{e^{M_{0}^{(4)} y_{0}^{4}}
-1}{\sigma^4 y_{0}^{4}} \biggr) \cdot\biggl( 1 + \frac{\kappa_4}{4!
\sigma^4} \biggr) + \max_{|x| \leq\kappa_4 y_0/4!} \biggl[\biggl| \frac
{F(x)}{x^2} \biggr| \biggr] \exp\bigl\{ M_{0}^{(4)} y_{0}^{4}\bigr\} \\
&&{} + \exp\bigl\{ M_{0}^{(4)} y_{0}^{4}\bigr\} \cdot\exp\biggl\{ \frac{\kappa_4
y_{0}^{4}}{4!}\biggr\} 4 \sigma^{-4}\\
&&{} + \exp\bigl\{ M_{0}^{(4)} y_{0}^{4}\bigr\} \frac
{\kappa_{4}^{2}}{3! 4! \sigma^8} \biggl( \frac{\exp\{
{\kappa_4 y_{0}^{4}}/{4!} \} - 1}{{\kappa_4 y_{0}^{4}}/{4!}} \biggr)
\end{eqnarray*}
and combine (\ref{eq:pippo1}), (\ref{eq:pippo2}), (\ref
{eq:pippo3}) and (\ref{eq:pippo4}), after noting that $|\xi|^5 (1 +
\xi
^4) + |\xi|^9 + |\xi|^3 + |\xi|^7 \leq4|\xi|^3 (1 + \xi^6)$ holds for
every $\xi$. Finally, in order to have the same multiplicative constant
in the right-hand sides of (\ref{eq:restopunt4}), (\ref
{eq:restopunt4fast}) and (\ref{eq:restopunt4der}), replace
$C_{4}^{\ast
\ast}$ and $4 C_{4}^{\ast\ast\ast}$ with $C_{4}^{\ast} := \max\{
C_{4}^{\ast\ast}; 4 C_{4}^{\ast\ast\ast} \}$.
\end{pf*}
\begin{pf*}{Proof of Lemma \protect\ref{lm:approssimazioni}} In view of
the independence of the random variables $X_{j,n}$ and (\ref
{eq:logtaylor}), one gets
\[
\log\psi_n(\xi) = -\frac{1}{2} \xi^2 + \sum_{r=2}^{\chi} (-1)^r
\frac
{\tilde{\lambda}_{r,n}}{(2r)!} \xi^{2r} + R_{k+\delta}(\xi),
\]
where
\[
R_{k+\delta}(\xi) := \sum_{j=1}^{n} \frac{c_{j,n}^{k} \xi
^k}{\sigma^k}
\epsilon_{k + \delta} \biggl(\frac{c_{j,n} \xi}{\sigma} \biggr) ,
\]
whence
%
%
\begin{equation} \label{eq:rina}
\psi_n(\xi) = e^{- \xi^2/2} \cdot\exp\Biggl\{\sum_{r=2}^{\chi} (-1)^r
\frac{\tilde{\lambda}_{r,n}}{(2r)!} \xi^{2r} \Biggr\} \cdot\exp\{
R_{k+\delta}(\xi)\} .
\end{equation}
Now, consider the function $z \mapsto f_{\xi}(z) = \exp\{
g_{\xi}(z)\}$ with
\[
g_{\xi}(z) := \sum_{r=1}^{\chi-1} (-1)^{r+1} \frac{\tilde{\lambda
}_{r+1,n}}{(2r + 2)!} \xi^{2(r+1)} z^r
\]
and its Taylor polynomial of order $(\chi-1)$ at $z=0$, say
$p_{\chi-1}(z)$. Then, recall the Fa\`{a} di Bruno formula, that is,
\begin{eqnarray*}
&&\frac{\ud^{(\chi)}}{\ud t^{(\chi)}} \exp\{(y(t))\} \\
&&\qquad= \sum_{(\ast)} \frac{\chi!}{k_1! k_2! \cdots k_{\chi}!} \exp\{
(y(t))\}
\biggl( \frac{y^{(1)}(t)}{1!} \biggr)^{k_1} \biggl( \frac{y^{(2)}(t)}{2!}
\biggr)^{k_2} \cdots\biggl( \frac{y^{(\chi)}(t)}{\chi!} \biggr)^{k_{\chi}}
\end{eqnarray*}
with $(\ast)$ meaning that the sum is carried out over all
nonnegative integer solutions $(k_1, \ldots, k_{\chi})$ of the equation
$k_1 + 2k_2 + \cdots+ \chi k_{\chi} = \chi$. An application of this
formula entails
that
\[
p_{\chi-1}(z) = 1 + \sum_{r=1}^{\chi-1} \tilde{P}_{r,n}(\xi) z^r,
\]
the functions $\tilde{P}_{r,n}(\xi)$ having been defined in
(\ref{eq:P}). Thus, when $z=1$, the Lagrange remainder can be written
with a suitable $u \in[0,1]$ as
\[
\frac{1}{\chi!} f_{\xi}^{(\chi)}(u)
= f_{\xi}(u) \sum_{(\ast)} \frac{1}{k_1! k_2! \cdots k_{\chi}!} \biggl(
\frac{g_{\xi}^{(1)}(u)}{1!} \biggr)^{k_1} \biggl( \frac{g_{\xi
}^{(2)}(u)}{2!} \biggr)^{k_2} \cdots\biggl( \frac{g_{\xi}^{(\chi
)}(u)}{\chi!} \biggr)^{k_{\chi}},
\]
which, after repeated application of the multinomial formula,
leads to
\[
\biggl| \frac{1}{\chi!} f_{\xi}^{(\chi)}(u) \biggr| \leq f_{\xi}(u) \sum
_{(\ast)} \prod_{m=1}^{\chi-1}
\sum_{\{ l_1 + \cdots+ l_{\chi-m} = k_m\}} |A_{1,m}^{l_1}(\xi)
\cdots
A_{\chi-m,m}^{l_{\chi-m}}(\xi)|
\]
with $A_{h,m}(\xi) := (-1)^{h+m} {h+m-1\choose m} \frac{\tilde
{\lambda}_{h+m,n}}{(2(h + m))!} \xi^{2(h+m)}$. We can then introduce
the quantity
\[
W_{\chi} := \Biggl[ \prod_{s=2}^{\chi} \max\biggl\{\frac{\kappa_{2s}}{\sigma
^{2s}}; 1\biggr\} \Biggr]^{\chi}
\]
to obtain, after an application of the Lyapunov inequality,
\begin{eqnarray*}
&&\sum_{\{ l_1 + \cdots+ l_{\chi-m} = k_m\}} |A_{1,m}^{l_1}(\xi)
\cdots
A_{\chi-m,m}^{l_{\chi-m}}(\xi)| \\[-1pt]
&&\qquad\leq\chi^{\chi} W_{\chi} \xi^{2mk_m}(\xi^2 + \xi^{k-2})^{k_m}
\cdot
\Biggl( \sum_{j=1}^{n} c_{j,n}^{k+2} \Biggr)^{{2mk_m}/{k}} ,
\end{eqnarray*}
whence
\[
\biggl| \frac{1}{\chi!} f_{\xi}^{(\chi)}(u) \biggr| \leq f_{\xi}(u)
\cdot\chi^{\chi^2} W_{\chi}^{\chi-1} \xi^k [(\xi^2 + \xi
^{k-2})^2 +
(\xi^2 + \xi^{k-2})^\chi] \cdot\Biggl(\sum_{j=1}^{n} c_{j,n}^{k+2}
\Biggr)
\]
and, using the bound $|\xi| \leq A_{k,\delta,n}$,
\[
|g_{\xi}(u)| \leq\sum_{s=2}^{\chi} \kappa_{2s} y_{0}^{2s} :=
B_{\chi}.
\]
Then,
%
%
\begin{eqnarray} \label{eq:tosca}
&&|\psi_n(\xi) - \eta_{k,n}(\xi)| \hspace*{-12pt}\nonumber\\[-1pt]
&&\qquad\leq e^{- \xi^2/2} \bigl\{ [f_{\xi}(1) - p_{\chi-1}(1)]
+ \bigl[e^{R_{k+\delta}(\xi)} - 1\bigr] \bigr\} \hspace*{-12pt}\nonumber\\[-9pt]\\[-9pt]
&&\qquad\leq e^{- \xi^2/2} \Biggl[ e^{B_{\chi}} \chi^{\chi^2} W_{\chi}^{\chi
-1} \xi^k [(\xi^2 + \xi^{k-2})^2 + (\xi^2 + \xi^{k-2})^\chi]
\cdot
\Biggl(\sum_{j=1}^{n} c_{j,n}^{k+2} \Biggr) \hspace*{-12pt}\nonumber\\[-1pt]
&&\qquad\quad\hspace*{58.1pt}{} + \biggl(\frac{\exp\{M_{0}^{(k+\delta)} y_{0}^{k}\} - 1}{M_{0}^{(k+\delta
)} y_{0}^{k}} \biggr) \frac{2 \overline{\mathfrak{m}}_{k+\delta}}{k!
\sigma^{k+\delta}} \Biggl(\sum_{j=1}^{n} |c_{j,n}|^{k+\delta} \Biggr) |\xi
|^{k+\delta} \Biggr].\hspace*{-12pt}\nonumber
\end{eqnarray}
After observing that $\xi^k [(\xi^2 + \xi^{k-2})^2 + (\xi^2 +
\xi^{k-2})^\chi] \leq|\xi|^{k+\delta} (1+\xi^2)[2\xi^2 + 2\xi^{2k-6}
+ 2^{\chi}\xi^{k-2} + 2^{\chi}\xi^{\chi k - k -2}]$ for every $\xi$,
one can take $p_{0,k}$ in (\ref{eq:restopunt}) to be equal to $1 +
(1+\xi^2)[2\xi^2 + 2\xi^{2k-6} + 2^{\chi}\xi^{k-2} + 2^{\chi}\xi
^{\chi
k - k -2}]$.

As for $|\psi_{n}^{\prime} - \eta_{k,n}^{\prime}|$, note that the inequality
%
%
\begin{eqnarray}\label{eq:traviata}
&&|\psi_{n}^{\prime}(\xi) - \eta_{k,n}^{\prime}(\xi)| \nonumber\\[-1pt]
&&\qquad\leq|\xi| \cdot|\psi_{n}(\xi) - \eta_{k,n}(\xi)| + e^{- \xi
^2/2} \biggl|
\frac{\ud}{\ud\xi} f_{\xi}(1) \biggr| \cdot| \exp\{R_{k+\delta
}(\xi)\} - 1 | \nonumber\\[-8.5pt]\\[-8.5pt]
&&\qquad\quad{} + e^{- \xi^2/2} |f_{\xi}(1)| \exp\{R_{k+\delta}(\xi)\} |
R_{k+\delta}^{\prime}(\xi) |\nonumber\\[-1pt]
&&\qquad\quad{} + e^{- \xi^2/2} \biggl|
\frac{\ud}{\ud\xi} \bigl( f_{\xi}(1) - p_{\chi-1}(1) \bigr)
\biggr|\nonumber
\end{eqnarray}
obtains. As regards the first summand, it will suffice to
multiply the upper bound stated in (\ref{eq:tosca}) for $|\psi_{n} -
\eta_{k,n}|$ by $|\xi|$. The latter factor in the second addend of
(\ref
{eq:traviata}) can be dominated by the last addend in (\ref{eq:tosca}),
while, for the former factor, one has
\[
\biggl| \frac{\ud}{\ud\xi} f_{\xi}(1) \biggr| \leq\exp\{B_{\chi}\}
\sum_{r=1}^{\chi-1} \frac{\kappa_{2r+2} y_{0}^{2r+1}}{(2r+1)!
\sigma} .
\]
As for the third addend, recall that $|f_{\xi}(1)| \leq\exp\{
B_{\chi}\}$ and $|R_{k+\delta}(\xi)| \leq y_{0}^{k} M_{0}^{(k+\delta
)}$. Moreover, $|R_{k+\delta}^{\prime}(\xi)| \leq\sum_{j=1}^{n}
\sigma^{-k}
\xi^{k-1} |c_{j,n}|^{k} \{ k |\epsilon_{k+\delta}(c_{j,n}\xi/\sigma
)| +\break
|\xi\sigma^{-1} c_{j,n}\times \epsilon_{k+\delta}^{\prime} (c_{j,n}\xi
/\sigma
)|\}
$ and, in view of Theorem 1 in Section 8.4 of Chow and Teicher (\citeyear{CT97}),
$(|\epsilon_{k+\delta}(x)| + |x \epsilon_{k+\delta}^{\prime}(x)|)
\leq4
\overline{\mathfrak{m}}_{k+\delta} |x|^{\delta}/(k-1)!$. It remains to
deal with the last summand in (\ref{eq:traviata}). Since $\frac
{\partial
}{\partial\xi}p_{\chi-1}$ is a Taylor polynomial for $\frac
{\partial
}{\partial\xi}f_{\xi}$, one can use the Bernstein integral form of the
remainder to obtain
\begin{eqnarray*}
&&\hspace*{-3pt}\biggl| \frac{\partial}{\partial\xi} \bigl(f_{\xi}(1) - p_{\chi
-1}(1) \bigr) \biggr|\\
&&\hspace*{-3pt}\qquad \leq\frac{1}{(\chi-1)!} \int_{0}^{1}
(1-u)^{\chi-1} \sum_{l=0}^{\chi} \pmatrix{\chi\cr l} \biggl| \frac
{\partial
^{l}}{\partial u^l} f_{\xi}(u) \biggr| \,\ud u \\
&&\hspace*{-3pt}\qquad\quad{}\times\sum_{j=1}^{n} c_{j,n}^{2(\chi-l+1)} \sum_{r=\chi-l}^{\chi
-1} |\xi
|^{2r + 1} \frac{\kappa_{2(r+1)}}{\sigma^{2(r+1)}} \\
&&\hspace*{-3pt}\qquad\leq\sum_{l=0}^{\chi} \frac{1}{(\chi-l)!} e^{B_{\chi}} \sum
_{(\ast
)_l} \prod_{m=1}^{l} \frac{1}{k_m!}
\Biggl( \frac{1}{m!} \sum_{r=m}^{\chi-1} \frac{\kappa_{2(r+1)}}{\sigma
^{2(r+1)}} |\xi|^{2r+2} \Biggr)^{k_m} \\
&&\hspace*{-3pt}\qquad\quad\hspace*{82.6pt}{}\times\sum_{r=\chi-l}^{\chi-1} |\xi|^{2r + 1} \frac{\kappa
_{2(r+1)}}{\sigma^{2(r+1)}} \Biggl( \sum_{j=1}^{n} c_{j,n}^{2l+2}
\Biggr) \cdot\Biggl( \sum_{j=1}^{n} c_{j,n}^{2(\chi-l+1)} \Biggr).
\end{eqnarray*}
To conclude, think of the last two sums of the $c_{j,n}$'s as
moments of order $2l$ and $2(\chi-l)$, respectively, and apply the
Lyapunov inequality to each sum to write
\[
\Biggl( \sum_{j=1}^{n} c_{j,n}^{2l+2} \Biggr) \cdot\Biggl( \sum
_{j=1}^{n} c_{j,n}^{2(\chi-l+1)} \Biggr) \leq\sum_{j=1}^{n}
c_{j,n}^{k+2} .
\]
\upqed\end{pf*}

\subsection{A complement to the proof of Theorem \protect\ref{thm:invalicabilita}}\label{app3}
We clarify why\break  $\inf_{t \in[0, \hat{t}]} \,\ud
_{\mathrm{TV}}(\mu(\cdot, t); \gamma_{\sigma})$ must be strictly positive under
the hypothesis that $\kappa_4(\tilde{\mu}_0)$ is different from zero.
Suppose, on the contrary, that $\inf_{t \in[0, \hat{t}]} \,\ud
_{\mathrm{TV}}(\mu
(\cdot, t);\break \gamma_{\sigma}) = 0$. Then, as $t \mapsto\ud_{\mathrm{TV}}(\mu
(\cdot, t); \gamma_{\sigma})$ is continuous on $[0, +\infty[$, by the
Wild expansion, there exists $t^{\ast}$ in $[0, \hat{t}]$ such that
$\ud
_{\mathrm{TV}}(\mu(\cdot, t^{\ast}); \gamma_{\sigma}) = 0$. On the one
hand, if
$t^{\ast} = 0$, then $\mu_0$ coincides with $\gamma_{\sigma}$ and this
contradicts the hypothesis that $\kappa_4(\tilde{\mu}_0)$ is different
from zero. On the other hand, if $t^{\ast} > 0$, then one can conclude,
in view of the Wild expansion, that $\mu_0$ possesses moments of every
order and is symmetric. A direct consequence of (\ref{eq:kac})
is that $\mathfrak{m}_{2k}(t) := \int_{\mathbb{R}} x^{2k} \mu(\ud x,
t)$ satisfies an ordinary first order differential equation, which
admits the constant $\int_{\mathbb{R}} x^{2k} \gamma_{\sigma}(\ud x)$
as a stationary solution. Hence, since we are assuming that $\mathfrak
{m}_{2k}(t^{\ast})$ is equal to such a constant, the uniqueness of the
solutions of the equations under consideration implies that $\mathfrak
{m}_{2k}(t) = \int_{\mathbb{R}} x^{2k} \gamma_{\sigma}(\ud x)$ for
every $t$ in $[0, \infty[$ and every positive integer $k$. In other
words, $\mu_0$ coincides with $\gamma_{\sigma}$, which once again
contradicts the fact that $\kappa_4(\tilde{\mu}_0)$ is different
from zero.

\subsection{The proofs of (\protect\ref{eq:mattia}) and (\protect\ref
{eq:paolo})}\label{app4} The proofs of (\ref{eq:mattia}) and (\ref
{eq:paolo}) follow from the following proposition. Let $g\dvtx
\mathbb{R} \rightarrow[0, +\infty[$ be an integrable function
and $\epsilon\dvtx\mathbb{R} \rightarrow\mathbb{R}$ be a
continuous, bounded function with $\epsilon(0) = 0$. Then
\[
\lim_{t \rightarrow+\infty} H(t) := e^{({1/4})t} \mathsf{E}_t
\Biggl\{ \Biggl( \int_{\mathbb{R}} g(\xi) \Biggl[ \sum_{j=1}^{\nu} \pi
_{j}^{4} | \epsilon(\pi_j \xi) | \Biggr]^2 \,\ud\xi
\Biggr)^{1/2} \Biggr\} = 0 .
\]
\begin{pf*}{Proofs of (\protect\ref{eq:mattia}) and (\protect\ref
{eq:paolo})} We fix an arbitrary small positive $\delta$
and show that there exists a value $t_{\delta}$ for which $|H(t)| <
\delta$, for every $t > t_{\delta}$.
First, in view of the fact that $\epsilon(\cdot)$ is continuous and $\epsilon(0) = 0$,
there exists a strictly posityve number $\overline{x}$ such that the inequality
\[
|\epsilon(x)| \leq\frac{\delta}{3 \sqrt{\Vert g \Vert_{1}}}
\]
holds for every $x$ in $[-\overline{x}, \overline{x}]$ with
$\Vert g \Vert_{1}= \int_{\mathbb{R}} g(\xi)
\,\ud\xi$. Set $\overline{\pi} := \max_{1 \leq j \leq\nu} \pi_j$
and $B
:= \overline{x}/|\overline{\pi}|$. $B$ is well defined since, due to
(\ref{eq:sfera}), $\overline{\pi} \neq0$. Now,
\begin{eqnarray*}
&&\Biggl\{\int_{\mathbb{R}} g(\xi) \Biggl[ \sum_{j=1}^{\nu} \pi_{j}^{4}
| \epsilon(\pi_j \xi) | \Biggr]^2 \,\ud\xi\Biggr\}^{1/2}
\\
&&\qquad\leq\Biggl\{\int_{{\{|\xi| \leq B\}}} g(\xi) \Biggl[ \sum_{j=1}^{\nu}
\pi_{j}^{4} | \epsilon(\pi_j \xi) | \Biggr]^2 \,\ud\xi
\Biggr\}^{1/2}\\
&&\qquad\quad{} + \Biggl\{\int_{{\{|\xi| \geq B\}}} g(\xi) \Biggl[ \sum
_{j=1}^{\nu} \pi_{j}^{4} | \epsilon(\pi_j \xi) | \Biggr]^2
\,\ud\xi\Biggr\}^{1/2} .
\end{eqnarray*}
For the integral over the internal region, one can write
\[
\Biggl\{\int_{{\{|\xi| \leq B\}}} g(\xi) \Biggl[ \sum_{j=1}^{\nu} \pi
_{j}^{4} | \epsilon(\pi_j \xi) | \Biggr]^2 \,\ud\xi\Biggr\}
^{1/2} \leq\frac{\delta}{3 \sqrt{\Vert g \Vert_{1}}} \Biggl( \sum
_{j=1}^{\nu} \pi_{j}^{4} \Biggr) \sqrt{\Vert g \Vert_{1}}
\]
and, taking expectation,
\[
e^{{(1/4)}t} \mathsf{E}_t \Biggl\{ \Biggl( \int_{{\{|\xi| \leq B\}
}} g(\xi) \Biggl[ \sum_{j=1}^{\nu} \pi_{j}^{4} | \epsilon(\pi_j \xi
) | \Biggr]^2 \,\ud\xi\Biggr)^{1/2} \Biggr\} \leq\delta/3 ,
\]
after a standard application of (\ref{eq:formulaGR}). At this
point, we define $M$ to be the maximum of $|\epsilon|$ and determine a
positive value $\overline{s}$ such that
\[
\int_{{\{|\xi| \geq\overline{s} \}}} g(\xi) \,\ud\xi\leq\biggl(
\frac{\delta}{3 M} \biggr)^2 .
\]
Given $S := \{ \omega| |\overline{\pi}(\omega)| <
\overline{x}/\overline{s} \}$, we write
\begin{eqnarray*}
&&e^{{(1/4)}t} \mathsf{E}_t \Biggl\{ \Biggl( \int_{{\{|\xi| \geq B\}
}} g(\xi) \Biggl[ \sum_{j=1}^{\nu} \pi_{j}^{4} | \epsilon(\pi_j \xi
) | \Biggr]^2 \,\ud\xi\Biggr)^{1/2} \Biggr\} \\
&&\qquad= e^{{(1/4)}t} \mathsf{E}_t \Biggl\{ \Biggl[ \int_{{\{|\xi| \geq
B\}}} g(\xi) \Biggl[ \sum_{j=1}^{\nu} \pi_{j}^{4} | \epsilon(\pi_j
\xi) | \Biggr]^2 \,\ud\xi\Biggr]^{1/2} ; S \Biggr\}
\\
&&\qquad\quad{} + e^{{(1/4)}t} \mathsf{E}_t \Biggl\{ \Biggl[ \int_{{\{|\xi| \geq
B\}}} g(\xi) \Biggl[ \sum_{j=1}^{\nu} \pi_{j}^{4} | \epsilon(\pi_j
\xi) | \Biggr]^2 \,\ud\xi\Biggr]^{1/2} ; S^c \Biggr\}.
\end{eqnarray*}
One can observe that $B(\omega) > \overline{s}$ for $\omega$
in $S$. We then have
\begin{eqnarray*}
&&e^{{(1/4)}t} \mathsf{E}_t \Biggl\{ \Biggl[ \int_{{\{|\xi| \geq B\}
}} g(\xi) \Biggl[ \sum_{j=1}^{\nu} \pi_{j}^{4} | \epsilon(\pi_j \xi
) | \Biggr]^2 \,\ud\xi\Biggr]^{1/2} ; S \Biggr\} \\
&&\qquad\leq e^{{(1/4)}t} \biggl\{ \int_{{\{|\xi| \geq\overline{s} \}}}
g(\xi) \,\ud\xi\biggr\}^{1/2} M \mathsf{E}_t \Biggl[ \sum_{j=1}^{\nu}
\pi_{j}^{4} \Biggr] \leq\delta/3 .
\end{eqnarray*}
For the remaining term,
\[
e^{{(1/4)}t} \mathsf{E}_t \Biggl\{ \Biggl[ \int_{{\{|\xi| \geq B\}
}} g(\xi) \Biggl[ \sum_{j=1}^{\nu} \pi_{j}^{4} | \epsilon(\pi_j \xi
) | \Biggr]^2 \,\ud\xi\Biggr]^{1/2} ; S^c \Biggr\} \leq e^{{(1/4)}t} M \sqrt
{\Vert
g \Vert_{1}} \mathsf{P}_t(S^c) .
\]
An application of Markov's inequality with exponent 6 yields
an upper bound for the probability of $S^c$, that is,
\[
\mathsf{P}_t(S^c) \leq\mathsf{E}_t [|\overline{\pi}|^6] \cdot
\biggl(\frac{\overline{s}}{\overline{x}} \biggr)^6 \leq\mathsf{E}_t
\Biggl[ \sum_{j=1}^{\nu} \pi_{j}^{6} \Biggr] \cdot\biggl(\frac{\overline
{s}}{\overline{x}} \biggr)^6 \leq e^{-{(3/8)}t} \cdot\biggl(\frac
{\overline{s}}{\overline{x}} \biggr)^6 .
\]
Hence,
\[
e^{{(1/4)}t} M \sqrt{\Vert g \Vert_{1}} \mathsf{P}_t(S^c) \leq
e^{-{(1/8)}t} M \sqrt{\Vert g \Vert_{1}} \cdot\biggl(\frac{\overline
{s}}{\overline{x}} \biggr)^6 .
\]
Taking $t_{\delta} = \max\{ -8 \log[ (\delta
/3 ) \cdot(\overline{x}/\overline{s} )^6 \cdot M^{-1}
\Vert g \Vert_{1}^{-1/2} ] ; 1 \}$ makes the right-hand side
of the last inequality smaller than $\delta/3$ for every $t >
t_{\delta
}$. This completes the proof.
\end{pf*}
\end{appendix}


%
\printaddresses

\end{document}